\documentclass[12pt]{amsart}
\usepackage{}
\usepackage{amsmath}
\usepackage{amsfonts}
\usepackage{amssymb}
\usepackage[all]{xy}           %xypic macro for latex2.09

\usepackage{bbding}
\usepackage{txfonts}
\usepackage{amscd}
\usepackage{xspace}

\usepackage[shortlabels]{enumitem}
\usepackage{ifpdf}
\ifpdf
  \usepackage[colorlinks,final,backref=page,hyperindex]{hyperref}
\else
  \usepackage[colorlinks,final,backref=page,hyperindex,hypertex]{hyperref}
\fi
\usepackage{tikz}
\usepackage[active]{srcltx}

\makeatletter

\newtheorem{thm}{Theorem}[section]
\newtheorem{lem}[thm]{Lemma}
\newtheorem{cor}[thm]{Corollary}
\newtheorem{pro}[thm]{Proposition}
\theoremstyle{definition}
\newtheorem{defi}[thm]{Definition}
\newtheorem{ex}[thm]{Example}
\newtheorem{rmk}[thm]{Remark}

\newcommand {\emptycomment}[1]{}

 \newcommand{\nc}{\newcommand}
 \nc{\delete}[1]{{}}
\nc{\CV}{\mathbf{C}}

\nc{\oprn}{\theta}
\nc{\Oprn}{\Theta}

\newcommand{\lon }{\,\rightarrow\,}
\newcommand{\be }{\begin{equation}}
\newcommand{\ee }{\end{equation}}

\nc{\frakg}{\mathfrak{g}}
\newcommand{\dgla}{{\rm dgLa}}
\newcommand{\gla}{{\rm gLa}}
\newcommand{\sgla}{{\rm sgLa}}
\newcommand{\g}{\mathfrak g}
\newcommand{\h}{\mathfrak h}

\nc{\adj}{\xspace}
\nc{\opt}{operator\xspace}
\nc{\Opt}{Operator\xspace}

\nc{\calo}{\mathcal{O}}
\newcommand{\huaP}{\mathcal{P}}
\newcommand{\huaC}{{\mathcal{C}}}%{\mathcal{C}}

\newcommand{\huaO}{{\mathcal{O}}}

\newcommand{\frkh}{\mathfrak h}

\newcommand{\half}{\frac{1}{2}}

\newcommand{\Courant}[1]{\left\llbracket  #1\right\rrbracket }

\newcommand{\Id}{\rm{Id}}

\newcommand{\br}[1]{   [ \cdot,    \cdot  ]   }

\newcommand{\dM}{\mathrm{d}}

\newcommand{\Hom}{\mathrm{Hom}}

\newcommand{\Der}{\mathsf{Der}}

\newcommand{\gl}{\mathfrak {gl}}
\newcommand{\sln}{\mathfrak {sl}}

\newcommand{\ad}{\mathsf{ad}}

\newcommand{\Img}{\mathrm{Im}}

\newcommand{\ZZ}{{\mathbb Z}}

\newcommand{\CC}{\mathbb{C}}

\begin{document}

\title[Homotopy $\calo$-operators and homotopy post-Lie algebras]{Homotopy Rota-Baxter operators, homotopy $\huaO$-operators  and homotopy post-Lie algebras}

\author{Rong Tang}
\address{Department of Mathematics, Jilin University, Changchun 130012, Jilin, China}
\email{tangrong16@mails.jlu.edu.cn}

\author{Chengming Bai}
\address{Chern Institute of Mathematics and LPMC, Nankai University,
Tianjin 300071, China}
\email{baicm@nankai.edu.cn}

\author{Li Guo}
\address{
Department of Mathematics and Computer Science,
         Rutgers University,
         Newark, NJ 07102}
\email{liguo@rutgers.edu}

\author{Yunhe Sheng}
\address{Department of Mathematics, Jilin University, Changchun 130012, Jilin, China}
\email{shengyh@jlu.edu.cn}

\vspace{-5mm}

\date{\today}

\begin{abstract}
Rota-Baxter operators, $\calo$-operators on Lie algebras and their interconnected pre-Lie and post-Lie algebras are important algebraic structures with applications in mathematical physics. This paper introduces the notions of a homotopy Rota-Baxter operator and a homotopy $\calo$-operator on a symmetric graded Lie algebra. Their characterization by Maurer-Cartan elements of suitable differential graded Lie algebras is provided. Through the action of a homotopy $\calo$-operator on a symmetric graded Lie algebra, we arrive at the notion of an \opt homotopy post-Lie algebra, together with its characterization in terms of Maurer-Cartan elements.
A cohomology theory of post-Lie algebras is established, with an application to 2-term skeletal operator homotopy post-Lie algebras.
\end{abstract}



\maketitle
\vspace{-10mm}
\tableofcontents

\allowdisplaybreaks

\vspace{-1.5cm}

\section{Introduction}
This paper initiates the homotopy study of Rota-Baxter operators, $\calo$-operators and the related pre-Lie algebras and post-Lie algebras.

\subsection{Background and motivation}

%\subsubsection{History of homotopy operations}

Homotopy is a fundamental notion in topology describing continuously deforming one function to another.

The first homotopy construction in algebra is the $A_\infty$-algebra of Stasheff, arising from his work on homotopy
characterization of connected based loop spaces~\cite{Sta63}. Later related developments include the work of Boardman and Vogt~\cite{BV68} about $E_\infty$-spaces on the infinite loop space,  the work of Schlessinger and Stasheff~\cite{SS85} about $L_\infty$-algebras on perturbations of rational homotopy types and the work of Chapoton and Livernet~\cite{CL} about pre-Lie$_\infty$ algebras, as well as homotopy Leibniz algebras~\cite{AP}. See~\cite{LV,MSS,Mer} for other examples.

Homotopy and operads are intimately related. In fact, operads were introduced as a tool in homotopy theory, specifically for iterated loop spaces.
Vaguely speaking, the homotopy of an algebraic structure is obtained when the defining relations of the algebraic structure is relaxed to hold up to homotopy. The resulting algebraic structure is homotopy equivalent to the original algebraic structure via a Homotopy Transfer Theorem. The idea
is behind much of the operadic developments. Through the work of Ginzburg-Kapranov~\cite{GK94},  Getzler-Jones~\cite{GJ94} and Markl~\cite{Ma1,MSS}, the homotopy of an operad $\huaP$ is in general defined to be the minimal model of $\huaP$. More precisely, $\huaP_\infty$ is the Koszul resolution as the cobar construction $\Omega \huaP^{\mathrm{i}}$ of the Koszul dual cooperad of $\huaP$~\cite{LV,MSS}. Since the latter notion makes sense only when $\huaP$ is quadratic, this approach does not apply to some important algebraic structures, such as the operad of Rota-Baxter algebras.

Rota-Baxter associative algebras were introduced in the probability study of G. Baxter and later found important applications in the Connes-Kreimer's algebraic approach to renormalization of quantum field theory \cite{CK}, among others. A Rota-Baxter operator on a Lie algebra is naturally the operator form of a classical $r$-matrix~\cite{STS} under certain conditions.
To better understand such connection in general, the notion of an $\calo$-operator (also called
a relative Rota-Baxter operator \cite{PBG} or a generalized Rota-Baxter operator \cite{Uch,Uch1})
on a Lie algebra was introduced by Kupershmidt~\cite{Ku},
which can be traced back to Bordemann
\cite{Bor}.

Both operators have been studied extensively in recent years~\cite{Gub}. Their operadic theory are challenging to establish since the operads are not quadratic and have nontrivial unary operations. At the same time, they provide promising testing grounds to expand the existing operad theory.

\subsection{Approach of the paper}

Because of the limitation of the Koszul resolution method to study homotopy of operads, other methods to give the related resolutions have been introduced.

Dotsenko and Khoroshkin~\cite{DK} used shuffle operads and the Gr\"obner bases method to show that, for the operad $ncRB$ of Rota-Baxter operators on associative algebras, the minimal model $ncRB_\infty$ is a quasi-free operad. They are able to write down formulas for small arities for differentials in the free resolutions for Quillen homology computation, though ``in general compact formulas
are yet to be found" as noted in the paper. One can expect a similarly challenging situation for the operad of Rota-Baxter operators on Lie algebras.

This paper follows the more basic and direct approach to homotopy via differential graded Lie algebras and Maurer-Cartan elements. This is in fact the approach taken in the early developments of characterizing algebraic and homotopy algebraic structures  before they are put under the more uniform and sophisticated framework of operads. These developments started with the well-known series of work by Gerstenhaber~\cite{Ge0,Ge} on deformations of associative algebras and by Nijenhuis and Richardson on Lie algebras~\cite{NR} a few years later.

This approach is based on the principle that objects of a certain algebraic structure on a vector space are given by degree 1 solutions of the Maurer-Cartan equation in a suitable differential graded Lie algebra built from the vector space. When the vector space is replaced by a graded vector space, similar solutions give objects in the homotopy algebraic structure. To make the idea more transparent, we recall the case of Lie algebras and homotopy Lie algebras (that is, $L_\infty$-algebras). Let $V$ be a vector space. Define the graded vector space
$\oplus_{n=0}^\infty \Hom(\wedge^nV,V)$
with the degree of elements in $\Hom(\wedge^nV,V)$ being $n-1$. For $f\in \Hom(\wedge^mV,V), g\in \Hom(\wedge^nV,V)$, define
$$ [f,g]_{NR}:=f\circ g- (-1)^{(m-1)(n-1)}g\circ f,$$
with $f\circ g\in \Hom(\wedge^{m+n-1}V,V)$ being defined by
\begin{equation}
(f\circ g)(v_1,\cdots,v_{m+n-1}):=\sum_{\sigma\in \mathbb S_{(n,m-1)}} (-1)^\sigma f(g(v_{\sigma(1)},\cdots,v_{\sigma(n)}),v_{\sigma(n+1)}, \cdots,v_{\sigma(m+n-1)}),
\label{eq:fgcirc}
\end{equation}
where the sum is over $(n,m-1)$-shuffles. Then $\big(\oplus_{n=0}^\infty \Hom(\wedge^nV,V),[\cdot,\cdot]_{NR}\big)$ is a differential graded Lie algebra with the trivial derivation. With this setup, a Lie algebra structure on $V$ is precisely a degree 1 solution $\omega\in \Hom(\wedge^2V,V)$ of the Maurer-Cartan equation
$$ [\omega,\omega ]_{NR}=0.$$
In the spirit of this characterization of Lie algebra structures by solutions of the Maurer-Cartan equation, homotopy Lie algebra structures can be characterized briefly as follows, with further details provided in Section~\ref{ss:prep}.
Given a graded vector space $V=\oplus_{i\in \ZZ} V^i$, let $S(V)$ denote the  symmetric algebra of $V$ and let $\Hom^n(S(V),V)$ denote the space of degree $n$ linear maps.
Define the Nijenhuis-Richardson bracket $[\cdot,\cdot]_{NR}$ on the graded vector space $\sum_{n\in \ZZ}\Hom^n(S(V),V)$ in a graded form of Eq.~\eqref{eq:fgcirc} (see Eq.~\eqref{eq:gfgcirc}). We again have a graded Lie algebra and (curved) $L_\infty$-algebras can be characterized as degree 1 solutions of the corresponding Maurer-Cartan equation.

This approach usually agrees with the operadic approach when both approaches apply. The more elementary nature of this approach makes it possible to be applied where the operadic approach cannot be applied yet, in particular to Rota-Baxter operators and its relative generalization, the $\huaO$-operators. Indeed, in~\cite{TBGS}, we took this approach to give a Maurer-Cartan characterization of $\calo$-operators (of weight zero) and further to establish a deformation theory and its controlling cohomology for $\huaO$-operators.

To illustrate our approach in a broader context, we focus on the Rota-Baxter Lie algebra for now and regard its algebraic structure (operad) as a pair $(\ell, T)$ consisting of a Lie bracket $\ell=[\cdot,\cdot]$ and a Rota-Baxter operator $T$. Then to obtain the homotopy of the Rota-Baxter Lie algebra, one can begin with taking homotopy of the binary operation $\ell$ or taking homotopy of the unary operation $T$. The homotopy of the Lie algebra $\ell_\infty=\{\ell_i\}_{i= 1}^{\infty}$ is well-known as the $L_\infty$-algebra~\cite{MSS}. Together with the natural Rota-Baxter operator action as defined in~\cite{PBG}, we have the Rota-Baxter homotopy Lie algebra $(\ell_\infty, T)$ or Rota-Baxter $L_\infty$-algebra. A Rota-Baxter Lie algebra naturally induces a post-Lie algebra, originated from an operadic study~\cite{Val} and found applications in mathematical physics and numerical analysis~\cite{BGN,Munthe-Kaas-Lundervold}. Likewise, a Rota-Baxter homotopy Lie algebra is expected to induce a homotopy post-Lie algebra whose construction is still not known beyond its conceptual definition as an operadic minimal model mentioned above, giving rise to the commutative diagram
\begin{equation}
\xymatrix{ (\ell, T) \ar[rr]\ar^{\text{RB action}}[d] && (\ell_\infty, T) \ar^{\text{RB action}}[d] \\
\text{post-Lie} \ar[rr]&& \text{homotopy post-Lie}
}
\label{eq:diaga}
\end{equation}
where the horizontal arrows are taking homotopy and the vertical arrows are taking actions of the Rota-Baxter operators.

In this paper, we will pursue the other direction, by taking homotopy of the Rota-Baxter operator $T$ and obtain $T_\infty:=\{T_i\}_{i= 0}^{\infty}$, without taking homotopy of $\ell$. We call the resulting structure $(\ell, T_\infty)$ the {\bf \opt homotopy Rota-Baxter Lie algebra} to distinguish it from the above mentioned Rota-Baxter homotopy Lie algebra. The action of the \opt homotopy Rota-Baxter operator $T_\infty$ gives rise to a variation of the homotopy post-Lie algebra which we will call the {\bf \opt homotopy post-Lie algebra} (see Remark~\ref{rk:postlie}). This gives another commutative diagram shown as the front rectangle in Eq.~\eqref{eq:diagb} while the diagram in Eq.~\eqref{eq:diaga} is embedded as the right rectangle.

Eventually, the {\em full} homotopy of the Rota-Baxter Lie algebra should come from the combined homotopies of both the Lie algebra structure and the Rota-Baxter operator structure, tentatively denoted by $(\ell_\infty, T_\infty)$ and called the {\bf full homotopy Rota-Baxter Lie algebra}. A suitable action of $T_\infty$ on $\ell_\infty$ should give the {\bf full homotopy post-Lie algebra} whose structure is still mysterious.
These various homotopies of the Rota-Baxter Lie algebra, as well as their derived homotopies of the post-Lie algebra, could be put together and form the following diagram
\begin{equation}
\xymatrix@!0{
&&& (\ell_\infty,T_\infty) \ar@{-}[rrrrrr]\ar@{-}'[d][dd]
&&& &&& (\ell_\infty, T) \ar@{-}[dd]
\\
(\ell,T_\infty) \ar@{-}[urrr]\ar@{-}[rrrrrr]\ar@{-}[dd]
&&&&& & (\ell,T) \ar@{-}[urrr]\ar@{-}[dd]&&&
\\
&&& \text{full homotopy} \atop \text{post-Lie} \ar@{-}'[rrr][rrrrrr]
&&& &&& \text{homotopy} \atop \text{post-Lie}
\\
\text{\opt homotopy} \atop \text{post-Lie} \ar@{-}[rrrrrr]\ar@{-}[urrr]
&&&&&& \text{post-Lie} \ar@{-}[urrr] &&&
}
\label{eq:diagb}
\end{equation}
where going in the left and the inside directions should be taking various homotopy, and going downward should be taking the actions of (homotopy) Rota-Baxter operators.

We note that when the weight of the Rota-Baxter operator is zero, the post-Lie algebra in the lower half of the diagram becomes a pre-Lie algebra and the homotopy post-Lie in the back-right-lower corner is a homotopy pre-Lie algebra. We find it quite amazing that the \opt homotopy post-Lie algebra in the front-left-lower corner is also the homotopy pre-Lie algebra (Corollary~\ref{co:oprelie}). This of course does not mean that the algebraic structure in the back-left-lower corner is also the homotopy pre-Lie algebra. It would be interesting to determine this structure even in this special case.

\subsection{Outline of the paper}

After some background on differential graded Lie algebras and Maurer-Cartan equations summarized in Section~\ref{ss:prep}, we introduce in Section~\ref{ss:hoop} the notion of a homotopy $\calo$-operator with any weight (Definition~\ref{de:homoop}). Using derived bracket, we construct a differential graded Lie algebra which characterizes homotopy $\calo$-operators of any weight as their Maurer-Cartan elements  (Theorem~\ref{hmotopy-o-operator-dgla}).

In Section~\ref{sec:post}, by applying a homotopy $\calo$-operator to a symmetric graded Lie algebra we obtain a variation of the homotopy post-Lie algebra, called the \opt homotopy post-Lie algebra. From here we can specialize in several directions and obtain interesting applications. First when the weight of the $\calo$-operator is taken to be zero, we obtain homotopy $\calo$-operators of weight zero. Since $\calo$-operators of weight zero naturally derive pre-Lie algebras~\cite{PBG}, it is expected that homotopy $\calo$-operators of weight zero derive homotopy pre-Lie algebras, that is pre-Lie$_\infty$-algebras. We confirm this in Corollary~\ref{co:oprelie}, yielding the commutative diagram
\begin{equation}
\xymatrix{
\calo\text{-operators} \ar^{\text{homotopy}}[rr] \ar_{\text{\adj}}[d] && \text{homotopy }\calo\text{-operators} \ar^{\text{\adj}}[d]\\
\text{pre-Lie} \ar^{\text{homotopy}}[rr] && \text{pre-Lie}_\infty
}
\label{eq:oprelie}
\end{equation}
In other words, the compositions of taking homotopy and taking operator action in either order gives the pre-Lie$_\infty$ algebras.

There has been quite much interest on constructions of post-Lie algebras in the recent literature~\cite{Burde16,Burde19,GLBJ,PLBG}. Another useful application of our general construction is the characterization of post-Lie algebra structures on a given Lie algebra using Maurer-Cartan elements in a suitable differential graded Lie algebra (Corollary~\ref{cor:dgla-post-lie}).

In Section~\ref{sec:repcoh} we first consider the cohomology theory of post-Lie algebras. In the ``abelian" case of pre-Lie algebras, the cohomology groups were first defined in~\cite{DA} by derived functors and then in~\cite{Mil} by resolutions of algebras from the Koszul duality theory in the framework of operads. An explicit cohomology theory of post-Lie algebras is not yet known. In this section we establish such a theory which reduces to the existing cohomology theory of pre-Lie algebras.
The third cohomology group of a post-Lie algebra are applied in Section~\ref{ss:deform} to classify 2-term skeletal \opt homotopy post-Lie algebras.

\medspace

\noindent
{\bf Notation.} We assume that all the vector spaces are over a field of characteristic zero. For a homogeneous element $x$ in a $\ZZ$-graded vector space, we also use $x$ in the exponent, as in $(-1)^x$, to denote its degree in order to simplify the notation.

\section{Homotopy $\huaO$-operators of weight $\lambda$}
\label{sec:hmtp}
In this section, we introduce the notion of a homotopy $\huaO$-operator of weight $\lambda$, where $\lambda$ is a constant. We construct a differential graded Lie algebra  (\dgla~for short) and show that  homotopy $\huaO$-operators of weight $\lambda$   can be characterized as its Maurer-Cartan elements to justify our definition.

\subsection{Maurer-Cartan elements and Nijenhuis-Richardson brackets}
\label{ss:prep}
We first recall some background needed in later sections.
\begin{defi}\label{MCE}{\rm (\cite{LV})}
  Let $(\g=\oplus_{k\in\mathbb Z}\g^k,[\cdot,\cdot],d)$ be a \dgla.  A degree $1$ element $\theta\in\g^1$ is called a {\bf Maurer-Cartan element} of $\g$ if it
  satisfies the following {\bf Maurer-Cartan equation}:
  \begin{equation}
  d \theta+\half[\theta,\theta]=0.
  \label{eq:mce}
  \end{equation}
  \end{defi}

A permutation $\sigma\in\mathbb S_n$ is called an $(i,n-i)$-shuffle if $\sigma(1)<\cdots <\sigma(i)$ and $\sigma(i+1)<\cdots <\sigma(n)$. If $i=0$ or $n$, we assume $\sigma=\Id$. The set of all $(i,n-i)$-shuffles will be denoted by $\mathbb S_{(i,n-i)}$. The notion of an $(i_1,\cdots,i_k)$-shuffle and the set $\mathbb S_{(i_1,\cdots,i_k)}$ are defined analogously.

Let $V=\oplus_{k\in\mathbb Z}V^k$ be a $\mathbb Z$-graded vector space. We will denote by $S(V)$ the symmetric  algebra  of $V$:
\begin{eqnarray*}
S(V):=\oplus_{i=0}^{\infty}S^i (V).
\end{eqnarray*}
Denote the product of homogeneous elements $v_1,\cdots,v_n\in V$ in $S^n(V)$ by $v_1\odot\cdots\odot v_n$. The degree of $v_1\odot\cdots\odot v_n$ is by definition the sum of the degree of $v_i$. For a permutation $\sigma\in\mathbb S_n$ and $v_1,\cdots, v_n\in V$,  the Koszul sign $\varepsilon(\sigma;v_1,\cdots,v_n)\in\{-1,1\}$ is defined by
\begin{eqnarray*}
v_1\odot\cdots\odot v_n=\varepsilon(\sigma;v_1,\cdots,v_n)v_{\sigma(1)}\odot\cdots\odot v_{\sigma(n)}
\end{eqnarray*}
and  the antisymmetric Koszul sign $\chi(\sigma;v_1,\cdots,v_n)\in\{-1,1\}$ is defined by
\begin{eqnarray*}
\chi(\sigma;v_1,\cdots,v_n):=(-1)^\sigma\varepsilon(\sigma;v_1,\cdots,v_n).
\end{eqnarray*}
Denote by $\Hom^n(S(V),V)$ the space of degree $n$ linear maps from the graded vector space $S(V)$ to the graded vector space $V$. Obviously, an element $f\in\Hom^n(S(V),V)$ is the sum of $f_i:S^i(V)\lon V$. We will write  $f=\sum_{i=0}^\infty f_i$.
 Set $C^n(V,V):=\Hom^n(S(V),V)$ and
$
C^*(V,V):=\oplus_{n\in\mathbb Z}C^n(V,V).
$
As the graded version of the classical Nijenhuis-Richardson bracket given in \cite{NR,NR2}, the {\bf Nijenhuis-Richardson bracket} $[\cdot,\cdot]_{NR}$ on the graded vector space $C^*(V,V)$ is given
by:
\begin{eqnarray}
[f,g]_{NR}:=f\circ g-(-1)^{mn}g\circ f,\,\,\,\,\forall f=\sum_{i=0}^\infty f_i\in C^m(V,V),~g=\sum_{j=0}^\infty g_j\in C^n(V,V),
\label{eq:gfgcirc}
\end{eqnarray}
where $f\circ g\in C^{m+n}(V,V)$ is defined by
 \begin{eqnarray}\label{NR-circ}
f\circ g&=&\Big(\sum_{i=0}^{\infty}f_i\Big)\circ\Big(\sum_{j=0}^{\infty}g_j\Big):=\sum_{k=0}^{\infty}\Big(\sum_{i+j=k+1}f_i\circ g_j\Big),
\end{eqnarray}
while $f_i\circ g_j\in \Hom(S^{i+j-1}(V),V)$ is defined by
\begin{eqnarray*}
(f_i\circ g_j)(v_1,\cdots,v_{i+j-1})
:=\sum_{\sigma\in\mathbb S_{(j,i-1)}}\varepsilon(\sigma)f_i(g_j(v_{\sigma(1)},\cdots,v_{\sigma(j)}),v_{\sigma(j+1)},\cdots,v_{\sigma(i+j-1)}),
\end{eqnarray*}
with the convention that $f_0\circ g_j:=0$\footnote{The linear map $f_0$ is just a distinguished element $\Phi\in V^0$.} and
\begin{eqnarray*}
(f_j\circ g_0)(v_1,\cdots,v_{j-1}):=f_j(g_0, v_{1},\cdots,v_{j-1}).
\end{eqnarray*}

The notion of a curved $L_\infty$-algebra was introduced in \cite{KS,Ma2}. See also \cite{Laz} for more applications.

\begin{thm}\label{graded-Nijenhuis-Richardson-bracket}{\rm (\cite{ALN,Ma2})}
With the above notations, $(C^*(V,V),[\cdot,\cdot]_{NR})$ is a graded Lie algebra (\gla~ for short). Its Maurer-Cartan elements $\sum_{i=0}^{\infty}l_i$ are the curved $L_\infty$-algebra structures on $V$.
\end{thm}

 We denote a curved $L_\infty$-algebra by $(V,\{l_k\}_{k=0}^\infty)$. A curved $L_\infty$-algebra $(V,\{l_k\}_{k=0}^\infty)$ with $l_0=0$ is exactly an $L_\infty$-algebra~\cite{LM,LS,stasheff:shla}.

\begin{defi}{\rm (\cite{ALN})}\label{symmetric-lie}
A {\bf symmetric graded Lie algebra (sgLa)} is a $\mathbb Z$-graded vector space $\g$ equipped with a bilinear bracket $[\cdot,\cdot]_\g:\g\otimes \g\lon \g$  of degree $1$, satisfying
\begin{itemize}
\item[\rm(1)] {\bf (graded symmetry)} $[x,y]_\g=(-1)^{xy}[y,x]_\g,$
\item[\rm(2)] {\bf (graded Leibniz rule)} $[x,[y,z]_\g]_\g=(-1)^{x+1}[[x,y]_\g,z]_\g+(-1)^{(x+1)(y+1)}[y,[x,z]_\g]_\g$.
\end{itemize}
Here $x,y,z$ are homogeneous elements in $\g$, which also denote their degrees when in exponent.
\end{defi}

A \sgla~ $(\g,[\cdot,\cdot]_\g)$ is just a curved $L_\infty$-algebra $(\g,\{l_k\}_{k=0}^\infty)$, in which $l_k=0$ for all $k\ge 0$ except $k=2$.

We recall the notion of the suspension and desuspension operators. Let $V=\oplus_{i\in\mathbb Z} V^i$ be a graded vector space, we define the {\bf suspension operator} $s:V\mapsto sV$ by assigning  $V$ to the graded vector space $sV=\oplus_{i\in\mathbb Z}(sV)^i$ with $(sV)^i:=V^{i-1}$.
There is a natural degree $1$ map $s:V\lon sV$ that is the identity map of the underlying vector space, sending $v\in V$ to its suspended copy $sv\in sV$. Likewise, the {\bf desuspension operator} $s^{-1}$ changes the grading of $V$ according to the rule $(s^{-1}V)^i:=V^{i+1}$. The  degree $-1$ map $s^{-1}:V\lon s^{-1}V$ is defined in the obvious way.

\begin{ex}
Let $V$ be a graded vector space. Then $s^{-1}\gl(V)$ is a \sgla~ where the symmetric Lie bracket is given by
\begin{eqnarray}\label{representation}
[s^{-1}f,s^{-1}g]:=(-1)^ms^{-1}[f,g],\quad\forall f\in\Hom^m(V,V),~~g\in\Hom^n(V,V).
\end{eqnarray}
\end{ex}

Let $(\g,[\cdot,\cdot]_\g)$ and $(\g',[\cdot,\cdot]_{\g'})$ be \sgla's. A {\bf homomorphism} from $\g$ to $\g'$ is a linear map $\phi:\g\lon \g'$ of degree $0$ such that
\begin{eqnarray*}
\phi([x_1,x_2]_\g)=[\phi(x_1),\phi(x_2)]_{\g'},\quad\forall x_1,x_2\in \g.
\end{eqnarray*}

\begin{defi}
 A linear map of graded vector spaces $D:\mathfrak{g}\lon \mathfrak{g}$ of degree $n$ is called a {\bf derivation of degree $n$} on a \sgla~   $(\g,[\cdot,\cdot]_\g)$ if
 \begin{eqnarray*}
 D[x,y]_\g=(-1)^n[Dx,y]_\g+(-1)^{n(x+1)}[x,Dy]_\g,\quad\forall x,y\in\g.
 \end{eqnarray*}
\end{defi}
We denote the vector space of derivations of degree $n$ by $\Der^n(\g)$.  Denote by $\Der(\g)=\oplus_{n\in\mathbb Z}\Der^n(\g)$, which is a graded vector space.

\begin{rmk}
A derivation of degree $n$ on a \sgla~ $(\g,[\cdot,\cdot]_\g)$  is just a homotopy derivation $\{\theta_k\}_{k=1}^{\infty}$ of degree $n$ on $(\g,[\cdot,\cdot]_\g)$, in which $\theta_k=0$ for all $k\ge 2$. See \cite{DL,KS} for more details.
\end{rmk}

By a straightforward check, we obtain
\begin{pro}
With the above notations, $(s^{-1}\Der(\g),[\cdot,\cdot])$ is a symmetric graded Lie subalgebra of $(s^{-1}\gl(\g),[\cdot,\cdot])$, where the bracket $[\cdot,\cdot]$ is defined by \eqref{representation}.
\end{pro}

\begin{defi}
An {\bf action} of a \sgla~$(\g,[\cdot,\cdot]_\g)$ on a \sgla~ $(\h,[\cdot,\cdot]_\h)$ is a homomorphism of graded vector spaces $\rho:\g\lon \Der(\h)$ of degree $1$ such that $s^{-1}\circ\rho:\g\lon s^{-1}\Der(\h)$ is a \sgla~ homomorphism.
\end{defi}

In particular, if $(\h,[\cdot,\cdot]_\h)$ is abelian, we obtain an action of a \sgla~ on a graded vector space.
It is obvious that $\ad:\g\lon\Der(\frak g)$ is an action of the \sgla~ $(\g,[\cdot,\cdot]_\g)$ on  itself, which is called the {\bf adjoint action}.

Let $\rho$ be an action of a \sgla~ $(\g,[\cdot,\cdot]_\g)$ on a graded vector space $V$. For $x\in \g^i$, we have $\rho(x)\in \Hom^{i+1}(V,V)$. Moreover, there is a \sgla~ structure on the direct sum $\g\oplus V$ given by
\begin{eqnarray*}
[x_1+v_1,x_2+v_2]_{\rho}:=[x_1,x_2]_\g+\rho(x_1)v_2+(-1)^{x_1x_2}\rho(x_2)v_1,\quad\forall x_1,x_2\in\g,v_1,v_2\in V.
\end{eqnarray*}
This \sgla~ is called the
{\bf semidirect product} of the \sgla~ $(\g,[\cdot,\cdot]_\g)$  and $(V;\rho)$, and denoted by $\g\ltimes_{\rho}V$.

\subsection{Homotopy $\calo$-operators of weight $\lambda$}
\label{ss:hoop}

Now we are ready to give the main notion of this paper.

\begin{defi}
Let  $\rho$ be an action of a \sgla~ $(\g,[\cdot,\cdot]_\g)$ on a \sgla~ $(\h,[\cdot,\cdot]_\h)$. A degree $0$ element $T=\sum_{i=0}^{+\infty}T_i\in \Hom(S(\h),\g)$ with $T_i\in \Hom(S^i(\h),\g)$ is called  {\bf a homotopy $\huaO$-operator of weight $\lambda$} on a \sgla~ $(\g,[\cdot,\cdot]_\g)$ with respect to the action $\rho$ if the following equalities hold for all $p\geq 0$ and all homogeneous elements $v_1,\cdots,v_p\in V$,
\begin{eqnarray}
\nonumber&&\sum_{1\le i< j\le p}(-1)^{v_i(v_1+\cdots+v_{i-1})+v_j(v_1+\cdots+v_{j-1})+v_iv_j}T_{p-1}(\lambda[v_i,v_j]_{\h},v_1,\cdots,\hat{v}_i,\cdots,\hat{v}_{j},\cdots,v_{p})\\
\label{homotopy-rota-baxter-o}&&+\sum_{k+l=p+1}\sum_{\sigma\in \mathbb S_{(l,1,p-l-1)}}\varepsilon(\sigma)T_{k-1}\Big(\rho\big(T_l(v_{\sigma(1)},\cdots,v_{\sigma(l)})\big)v_{\sigma(l+1)},v_{\sigma(l+2)},\cdots,v_{\sigma(p)}\Big)\\
&=&\frac{1}{2}\sum_{k+l=p+1}\sum_{\sigma\in \mathbb S_{(k-1,l)}}\varepsilon(\sigma)[T_{k-1}(v_{\sigma(1)},\cdots,v_{\sigma(k-1)}),T_l(v_{\sigma(k)},\cdots,v_{\sigma(p)})]_\g.
%&=&\frac{1}{2}\sum_{k+l=s+1}\sum_{\sigma\in \mathbb S_{(s-l,l)}}\varepsilon(\sigma)[T_{k-1}(v_{\sigma(1)},\cdots,v_{\sigma(s-l)}),T_l(v_{\sigma(s-l+1)},\cdots,v_{\sigma(s)})]_\g.
\nonumber
\end{eqnarray}
\label{de:homoop}
\end{defi}

\begin{rmk}
The linear map $T_0$ is just a element $\Omega\in \g^0$. Below are the generalized Rota-Baxter identities for $p=0,1,2:$
\begin{eqnarray*}
\label{weak-1}&&[\Omega,\Omega]_\g=0,\\
\label{weak-2}&&T_1(\rho(\Omega)v_1)=[\Omega,T_1(v_1)]_\g,\\
\label{weak-3}&&[T_1(v_1),T_1(v_2)]_\g-T_1\Big(\rho(T_1(v_1))v_2+(-1)^{v_1v_2}\rho(T_1(v_2))v_1+\lambda[v_1,v_2]_\h\Big)\\
&&=T_2(\rho(\Omega)v_1,v_2)+(-1)^{v_1v_2}T_2(\rho(\Omega)v_2,v_1)-[\Omega,T_2(v_1,v_2)]_\g.
\end{eqnarray*}
\end{rmk}

\begin{rmk}
  If the \sgla~  reduces to a Lie algebra and the action reduces to an action of a Lie algebra on another Lie algebra, the above definition reduces to the definition of an $\huaO$-operator of weight $\lambda$ on a Lie algebra. More precisely, the linear map $T:\h \longrightarrow\g$ satisfies
$$   [Tu,Tv]=T\big(\rho(Tu)(v)-\rho(Tv)(u)+\lambda[u,v]_\h\big),\quad\forall u,v\in \h.
$$
\end{rmk}

\begin{defi}
A degree $0$ element $R=\sum_{i=0}^{+\infty}R_i\in \Hom(S(\g),\g)$ with $R_i\in \Hom(S^i(\g),\g)$ is called  {\bf a homotopy Rota-Baxter operator of weight $\lambda$} on a \sgla~ $(\g,[\cdot,\cdot]_\g)$ if the following equalities hold for all $p\geq 0$ and all homogeneous elements $x_1,\cdots,x_p\in \g$,
\begin{eqnarray*}
&&\sum_{1\le i< j\le p}(-1)^{x_i(x_1+\cdots+x_{i-1})+x_j(x_1+\cdots+x_{j-1})+x_ix_j}R_{p-1}(\lambda[x_i,x_j]_{\g},x_1,\cdots,\hat{x}_i,\cdots,\hat{x}_{j},\cdots,x_{p})\\
&&+\sum_{k+l=p+1}\sum_{\sigma\in \mathbb S_{(l,1,p-l-1)}}\varepsilon(\sigma)R_{k-1}\big([R_l(x_{\sigma(1)},\cdots,x_{\sigma(l)}),x_{\sigma(l+1)}]_\g,x_{\sigma(l+2)},\cdots,x_{\sigma(p)}\big)\\
&=&\frac{1}{2}\sum_{k+l=p+1}\sum_{\sigma\in \mathbb S_{(k-1,l)}}\varepsilon(\sigma)[R_{k-1}(x_{\sigma(1)},\cdots,x_{\sigma(k-1)}),R_l(x_{\sigma(k)},\cdots,x_{\sigma(p)})]_\g.
%&=&\frac{1}{2}\sum_{k+l=s+1}\sum_{\sigma\in \mathbb S_{(s-l,l)}}\varepsilon(\sigma)[R_{k-1}(x_{\sigma(1)},\cdots,x_{\sigma(s-l)}),R_l(x_{\sigma(s-l+1)},\cdots,x_{\sigma(s)})]_\g.
\end{eqnarray*}
\end{defi}

\begin{rmk}
A homotopy Rota-Baxter operator $R=\sum_{i=0}^{+\infty}R_i\in \Hom(S(\g),\g)$ of weight $\lambda$ on a \sgla~ $(\g,[\cdot,\cdot]_\g)$ is a homotopy $\huaO$-operator of weight $\lambda$ with respect to the adjoint action $\ad$.
If moreover the \sgla~  reduces to a Lie algebra, then the resulting linear operator $R:\g\longrightarrow \g$ is a {\bf Rota-Baxter operator of weight $\lambda$} in the sense that
$$ [R(x),R(y)]_\g=R\big([R(x),y]_\g+ [x,R(y)]_\g +\lambda  [x,y]_\g\big), \quad \forall x, y \in \g.
$$
\end{rmk}

In the sequel, we construct a \dgla~ and show that  homotopy $\huaO$-operators of weight $\lambda$ can be characterized as its Maurer-Cartan elements to justify our definition of homotopy $\huaO$-operators of weight $\lambda$.
For this purpose, we recall the derived bracket construction of graded Lie algebras. Let $(\g,[\cdot,\cdot]_\g,d)$ be a \dgla. We define a new bracket on $s\g$  by
\begin{eqnarray}
[sx,sy]_{d}:=(-1)^{x}s[dx,y]_\g,\quad \forall x,y\in\g.
\end{eqnarray}
The new bracket is called the  {\bf derived bracket} \cite{Kosmann-Schwarzbach,Vor}. It is well known that the derived bracket is a graded Leibniz bracket on the shifted graded space $s\g$. Note that the derived bracket is not graded skew-symmetric in general. We recall a basic result.

\begin{pro}{\rm (\cite{Kosmann-Schwarzbach})}\label{old-derived}
Let $(\g,[\cdot,\cdot]_\g,d)$ be a \dgla, and let $\h\subset \g$ be a subalgebra which is abelian, i.e. $[\h,\h]_\g=0$. If the derived bracket is closed on $s\h$, then $(s\h,[\cdot,\cdot]_d)$ is a \gla.
\end{pro}

Let $\rho$ be an action of a \sgla~ $(\g,[\cdot,\cdot]_\g)$ on a \sgla~ $(\h,[\cdot,\cdot]_\h)$. Consider the graded vector space
$C^*(\h,\g):=\oplus_{n\in\mathbb Z}\Hom^n(S(\h),\g)$.
Define a linear map $\dM:\Hom^n(S(\h),\g)\longrightarrow \Hom^{n+1}(S(\h),\g)$ by
\begin{eqnarray}
\nonumber&&(\dM g)_p ( v_1,\cdots, v_{p} )\\
\label{diff}&=&\sum_{1\le i< j\le p}(-1)^{n+1+v_i(v_1+\cdots+v_{i-1})+v_j(v_1+\cdots+v_{j-1})+v_iv_j} g_{p-1}(\lambda[v_i,v_j]_{\h},v_1,\cdots,\hat{v}_i,\cdots,\hat{v}_{j}, \cdots,v_{p}).
\end{eqnarray}
Also define a graded bracket operation
 $$\Courant{\cdot,\cdot}: \Hom^m(S(\h),\g)\times \Hom^n(S(\h),\g)\longrightarrow \Hom^{m+n+1}(S(\h),\g)$$
by
\begin{eqnarray}
\nonumber&&\Courant{f,g}_p( v_1,\cdots,v_{p})\\
\label{graded-Lie}&=&-\sum_{k+l=p+1}\sum_{\sigma\in \mathbb S_{(l,1,p-l-1)}}\varepsilon(\sigma)f_{k-1}\Big(\rho\big(g_l(v_{\sigma(1)},\cdots,v_{\sigma(l)})\big)v_{\sigma(l+1)},v_{\sigma(l+2)},\cdots,v_{\sigma(p)}\Big)\\
\nonumber&&+(-1)^{(m+1)(n+1)}\sum_{k+l=s+1}\sum_{\sigma\in \mathbb S_{(k-1,1,p-k)}}\varepsilon(\sigma)g_l\Big(\rho\big(f_{k-1}(v_{\sigma(1)},\cdots,v_{\sigma(k-1)})\big)v_{\sigma(k)}, v_{\sigma(k+1)},\cdots,v_{\sigma(p)}\Big)\\
\nonumber&&-\sum_{k+l=p+1}\sum_{\sigma\in \mathbb S_{(k-1,l)}}(-1)^{n(v_{\sigma(1)}+\cdots+v_{\sigma(k-1)})+m+1}\varepsilon(\sigma)[f_{k-1}(v_{\sigma(1)},\cdots,v_{\sigma(k-1)}),g_l(v_{\sigma(k)},\cdots,v_{\sigma(p)})]_\g
\end{eqnarray}
for all  $f=\sum_if_i\in\Hom^m(S(\h),\g)$, $g=\sum g_i\in \Hom^n(S(\h),\g)$  with  $f_i, g_i\in\Hom(S^i(\h),\g)$
 and $v_1,\cdots, v_{p} \in\h.$ Here we write $\dM g=\sum_i(\dM g)_i$ with $(\dM g)_i\in\Hom(S^i(\h),\g)$, and $\Courant{f,g}=\sum_i\Courant{f,g}_i$ with $\Courant{f,g}_i\in\Hom(S^i(\h),\g)$.

\begin{thm}\label{dgla-deforamtion-homotopy}
Let $\rho$ be an action of a \sgla~ $(\g,[\cdot,\cdot]_\g)$ on a \sgla~ $(\h,[\cdot,\cdot]_\h)$. Then $(sC^*(\h,\g),\Courant{\cdot,\cdot},\dM)$ is a \dgla.
\end{thm}

\begin{proof}
By Theorem \ref{graded-Nijenhuis-Richardson-bracket}, the graded Nijenhuis-Richardson bracket $[\cdot,\cdot]_{NR}$ associated to the direct sum vector space $\g\oplus \h$ gives rise to a \gla~ $(C^*(\g\oplus\h,\g\oplus\h),[\cdot,\cdot]_{NR})$. Obviously
$$
C^*(\h,\g)=\bigoplus_{n\in\mathbb Z}\Hom^n(S(\h),\g)
$$
is an abelian subalgebra. We denote the symmetric graded Lie brackets $[\cdot,\cdot]_\g$ and $[\cdot,\cdot]_\h$ by $\mu_\g$ and $\mu_\h$ respectively.  Since $\rho$ is an  action of the \sgla~ $(\g,[\cdot,\cdot]_\g)$, $\mu_\g+\rho$ is a semidirect product \sgla~ structure on $\g\oplus\h$. By Theorem \ref{graded-Nijenhuis-Richardson-bracket}, we deduce that $\mu_\g+\rho$ and $\lambda\mu_\h$ are Maurer-Cartan elements of the \gla~ $(C^*(\g\oplus\h,\g\oplus\h),[\cdot,\cdot]_{NR})$. Define a differential $d_{\mu_\g+\rho}$ on $(C^*(\g\oplus\h,\g\oplus\h),[\cdot,\cdot]_{NR})$ via
$$
  d_{\mu_\g+\rho}:=[\mu_\g+\rho,\cdot]_{NR}.
$$
Further, we define the derived bracket on the graded vector space $\oplus_{n\in\mathbb Z}\Hom^n(S(\h),\g)$ by
\begin{eqnarray}
\label{d-bracket}\Courant{f,g}:=(-1)^{m}[d_{\mu_\g+\rho}f,g]_{NR}=(-1)^{m}[[\mu_\g+\rho,f]_{NR},g]_{NR},
\end{eqnarray}
for all $f=\sum f_i\in\Hom^m(S(\h),\g), ~g=\sum_ig_i\in\Hom^n(S(\h),\g).$ Write
 $$[\mu_\g+\rho,f]_{NR}=\sum_{i=0}^{\infty}[\mu_\g+\rho,f]_{NR}^i \quad \mbox{with} \quad [\mu_\g+\rho,f]_{NR}^i\in\Hom(S^i(\h),\g).$$
  By \eqref{NR-circ}, for all $k\ge 2$,  $x_1,\cdots,x_k\in\g$ and $v_1\cdots,v_k\in\h$, we have
\begin{eqnarray*}
&&[\mu_\g+\rho,f]_{NR}^k\Big((x_1,v_1),\cdots,(x_k,v_k)\Big)\\
&=&\Big((\mu_\g+\rho)\circ f_{k-1}-(-1)^mf_{k-1}\circ (\mu_\g+\rho)\Big)\Big((x_1,v_1),\cdots,(x_k,v_k)\Big)\\
&=&\sum_{i=1}^{k}(-1)^\alpha(\mu_\g+\rho)\Big(f_{k-1}\big((x_1,v_1),\cdots,\widehat{(x_i,v_i)},\cdots,(x_k,v_k)\big),(x_i,v_i)\Big)\\
&&-(-1)^m\sum_{1\le i<j\le k}(-1)^\beta f_{k-1}\Big((\mu_\g+\rho)\big((x_i,v_i),(x_j,v_j)\big),(x_1,v_1),\cdots,\widehat{(x_i,v_i)},\cdots,\widehat{(x_j,v_j)},\cdots,(x_k,v_k)\Big)\\
&=&\sum_{i=1}^{k}(-1)^\alpha(\mu_\g+\rho)\Big(\big(f_{k-1}(v_1,\cdots,\hat{v}_i,\cdots,v_k),0\big),(x_i,v_i)\Big)\\
&&-(-1)^m\hspace{-.4cm}\sum_{1\le i<j\le k}(-1)^\beta f_{k-1}\Big(\big([x_i,x_j]_\g,\rho(x_i)v_j+(-1)^{v_iv_j}\rho(x_j)v_i\big),(x_1,v_1),\cdots,\widehat{(x_i,v_i)},\cdots,\widehat{(x_j,v_j)},\cdots,(x_k,v_k)\Big)\\
&=&\sum_{i=1}^{k}(-1)^\alpha\big([f_{k-1}(v_1,\cdots,\hat{v}_i,\cdots,v_k),x_i]_\g,\rho(f_{k-1}(v_1,\cdots,\hat{v}_i,\cdots,v_k))v_i\big)\\
&&-(-1)^m\sum_{1\le i<j\le k}(-1)^\beta \big(f_{k-1}(\rho(x_i)v_j+(-1)^{v_iv_j}\rho(x_j)v_i,v_1,\cdots,\hat{v}_i,\cdots,\hat{v}_j,\cdots,v_k),0\big),
\end{eqnarray*}
here $\alpha=v_i(v_{i+1}+\cdots+v_k)$ and $\beta=v_i(v_1+\cdots+v_{i-1})+v_j(v_1+\cdots+v_{j-1})+v_iv_j$. On the other hand, we have $[\mu_\g+\rho,f]_{NR}^0=0$ and $[\mu_\g+\rho,f]_{NR}^1(x_1,v_1)=\big([f_0,x_1]_\g,\rho(f_0)v_1\big)$.

Moreover, we obtain
\begin{eqnarray*}
&&[[\mu_\g+\rho,f]_{NR},g]_{NR}^p\Big((x_1,v_1),\cdots,(x_p,v_p)\Big)\\
&=&\Big(\sum_{k+l=p+1}[\mu_\g+\rho,f]_{NR}^k\circ g_l-(-1)^{(m+1)n}\sum_{k+l=p+1}g_l\circ [\mu_\g+\rho,f]_{NR}^k\Big)\Big((x_1,v_1),\cdots,(x_p,v_p)\Big).
\end{eqnarray*}
By straightforward computations, we have
\begin{eqnarray*}
&&([\mu_\g+\rho,f]_{NR}^k\circ g_l)\Big((x_1,v_1),\cdots,(x_p,v_p)\Big)\\
&=&\sum_{\sigma\in \mathbb S_{(l,p-l)}}\varepsilon(\sigma)[\mu_\g+\rho,f]_{NR}^k\Big(g_l\big((x_{\sigma(1)},v_{\sigma(1)}),\cdots,
(x_{\sigma(l)},v_{\sigma(l)})\big),(x_{\sigma(l+1)},v_{\sigma(l+1)}),\cdots,(x_{\sigma(p)},v_{\sigma(p)})\Big)\\
&=&\sum_{\sigma\in \mathbb S_{(l,p-l)}}\varepsilon(\sigma)[\mu_\g+\rho,f]_{NR}^k\Big(\big(g_l(v_{\sigma(1)},\cdots,v_{\sigma(l)}),0\big),(x_{\sigma(l+1)},
v_{\sigma(l+1)}),\cdots,(x_{\sigma(p)},v_{\sigma(p)})\Big)\\
&=&\sum_{\sigma\in \mathbb S_{(l,p-l)}}\varepsilon(\sigma)(-1)^{\bar{\alpha}}\big([f_{k-1}(v_{\sigma(l+1)},\cdots,v_{\sigma(p)}),g_l(v_{\sigma(1)},\cdots,v_{\sigma(l)})]_\g,0\big)\\
&&-(-1)^m\sum_{\sigma\in \mathbb S_{(l,p-l)}}\varepsilon(\sigma)\sum_{j=l+1}^{p}(-1)^{\bar{\beta}}\Big(f_{k-1}\big(\rho(g_l(v_{\sigma(1)},\cdots,
v_{\sigma(l)}))v_{\sigma(j)},v_{\sigma(l+1)},\cdots,\hat{v}_{\sigma(j)},\cdots,v_{\sigma(p)}\big),0\Big),
\end{eqnarray*}
where $\bar{\alpha}=(v_{\sigma(1)}+\cdots+v_{\sigma(l)}+n)(v_{\sigma(l+1)}+\cdots+v_{\sigma(p)})$ and $\bar{\beta}=v_{\sigma(j)}(v_{\sigma(l+1)}+\cdots+v_{\sigma(j-1)})$. For any $\sigma\in\mathbb S_{(l,p-l)}$, we define $\tau=\tau_\sigma\in\mathbb S_{(s-l,l)}$ by
\[
\tau(i)=\left\{
\begin{array}{ll}
\sigma(i+l), & 1\le i\le p-l;\\
\sigma(i-p+l), & p-l+1\le i\le s.
\end{array}
\right.
\]
Thus $\varepsilon(\tau;v_1,\cdots,v_p)=\varepsilon(\sigma;v_1,\cdots,v_p)(-1)^{(v_{\sigma(1)}+\cdots+v_{\sigma(l)})(v_{\sigma(l+1)}+\cdots+v_{\sigma(p)})}$.
In fact, the elements of $\mathbb S_{(l,p-l)}$ are in bijection with the elements of $\mathbb S_{(p-l,l)}$. Moreover, by $k+l=p+1$, we have
\begin{eqnarray*}
&&\sum_{\sigma\in \mathbb S_{(l,p-l)}}\varepsilon(\sigma)(-1)^{\bar{\alpha}}\big([f_{k-1}(v_{\sigma(l+1)},\cdots,v_{\sigma(p)}),g_l(v_{\sigma(1)},\cdots,v_{\sigma(l)})]_\g,0\big)\\
&=&\sum_{\tau\in \mathbb S_{(k-1,l)}}\varepsilon(\tau)(-1)^{n(v_{\tau(1)}+\cdots+v_{\tau(k-1)})}\big([f_{k-1}(v_{\tau(1)},\cdots,v_{\tau(k-1)}),g_l(v_{\tau(k)},\cdots,v_{\tau(p)})]_\g,0\big).
\end{eqnarray*}
For any $\sigma\in\mathbb S_{(l,p-l)}$ and $l+1\le j\le p$, we define $\tau=\tau_{\sigma,j}\in\mathbb S_{(l,1,p-l-1)}$ by
\[
\tau(i)=\left\{
\begin{array}{ll}
\sigma(i), & 1\le i\le l;\\
\sigma(j), &  i=l+1;\\
\sigma(i-1),   & l+2\le i\le j;\\
\sigma(i),& j+1\le i\le p.
\end{array}
\right.
\]
Thus we have $\varepsilon(\tau;v_1,\cdots,v_p)=\varepsilon(\sigma;v_1,\cdots,v_p)(-1)^{v_{\sigma(j)}(v_{\sigma(l+1)}+\cdots+v_{\sigma(j-1)})}$. Then
\begin{eqnarray*}
&&\sum_{\sigma\in \mathbb S_{(l,p-l)}}\varepsilon(\sigma)\sum_{j=l+1}^{p}(-1)^{\bar{\beta}}\Big(f_{k-1}\big(\rho(g_l(v_{\sigma(1)},\cdots,v_{\sigma(l)}))v_{\sigma(j)},
v_{\sigma(l+1)},\cdots,\hat{v}_{\sigma(j)},\cdots,v_{\sigma(p)}\big),0\Big)\\
&=&\sum_{\tau\in \mathbb S_{(l,1,p-l-1)}}\varepsilon(\tau)\Big(f_{k-1}\big(\rho(g_l(v_{\tau(1)},\cdots,v_{\tau(l)}))v_{\tau(l+1)},v_{\tau(l+2)},\cdots,v_{\tau(p)}\big),0\Big).
\end{eqnarray*}
Therefore, we obtain
\begin{eqnarray*}
&&([\mu_\g+\rho,f]_{NR}^k\circ g_l)\Big((x_1,v_1),\cdots,(x_p,v_p)\Big)\\
&=&\sum_{\sigma\in \mathbb S_{(k-1,l)}}\varepsilon(\sigma)(-1)^{n(v_{\sigma(1)}+\cdots+v_{\sigma(k-1)})}\big([f_{k-1}(v_{\sigma(1)},\cdots,v_{\sigma(k-1)}),
g_l(v_{\sigma(k)},\cdots,v_{\sigma(p)})]_\g,0\big)\\
&&-(-1)^m\sum_{\sigma\in \mathbb S_{(l,1,p-l-1)}}\varepsilon(\sigma)\Big(f_{k-1}\big(\rho(g_l(v_{\sigma(1)},\cdots,v_{\sigma(l)}))v_{\sigma(l+1)},
v_{\sigma(l+2)},\cdots,v_{\sigma(p)}\big),0\Big).
\end{eqnarray*}

On the other hand,
\begin{eqnarray*}
&&\Big(g_l\circ [\mu_\g+\rho,f]_{NR}^k\Big)\Big((x_1,v_1),\cdots,(x_p,v_p)\Big)\\
&=&\sum_{\sigma\in \mathbb S_{(k,n-k)}}\varepsilon(\sigma)g_l\Big([\mu_\g+\rho,f]_{NR}^k\big((x_{\sigma(1)},v_{\sigma(1)}),\cdots,(x_{\sigma(k)},v_{\sigma(k)})\big),
(x_{\sigma(k+1)},v_{\sigma(k+1)}),\cdots,(x_{\sigma(p)},v_{\sigma(p)})\Big)\\
&=&\sum_{\sigma\in \mathbb S_{(k,p-k)}}\varepsilon(\sigma)\sum_{i=1}^{k}(-1)^{\alpha'}\Big(g_l\big(\rho(f_{k-1}(v_{\sigma(1)},\cdots,\hat{v}_{\sigma(i)},\cdots,
v_{\sigma(k)}))v_{\sigma(i)},v_{\sigma(k+1)}\cdots,v_{\sigma(p)}\big),0\Big),
\end{eqnarray*}
where $\alpha'=v_{\sigma(i)}(v_{\sigma(i+1)}+\cdots+v_{\sigma(k)})$. For any $\sigma\in\mathbb S_{(k,p-k)}$ and $1\le i\le k$, we define $\tau=\tau_{\sigma,i}\in\mathbb S_{(k-1,1,p-k)}$ by
\[
\tau(j)=\left\{
\begin{array}{ll}
\sigma(j), & 1\le j\le i-1;\\
\sigma(j+1), &  i\le j\le k-1;\\
\sigma(i),   & j=k;\\
\sigma(j),& k+1\le j\le p.
\end{array}
\right.
\]
Thus we have $\varepsilon(\tau;v_1,\cdots,v_p)=\varepsilon(\sigma;v_1,\cdots,v_p)(-1)^{v_{\sigma(i)}(v_{\sigma(i+1)}+\cdots+v_{\sigma(k)})}$. Then we have
\begin{eqnarray*}
&&\Big(g_l\circ [\mu_\g+\rho,f]_{NR}^k\Big)\Big((x_1,v_1),\cdots,(x_p,v_p)\Big)\\
&=&\sum_{\sigma\in \mathbb S_{(k-1,1,p-k)}}\varepsilon(\sigma)\Big(g_l\big(\rho(f_{k-1}(v_{\sigma(1)},\cdots,v_{\sigma(k-1)}))v_{\sigma(k)},v_{\sigma(k+1)}\cdots,v_{\sigma(p)}\big),0\Big).
\end{eqnarray*}
By \eqref{d-bracket}, we obtain that the derived bracket $\Courant{\cdot,\cdot}$ is closed on $sC^*(\h,\g)$, and given by \eqref{graded-Lie}. Therefore,   $(sC^*(\h,\g),\Courant{\cdot,\cdot})$ is a \gla.

Moreover, by $\Img\rho\subset\Der(\h)$, we have $[\mu_\g+\rho,\lambda\mu_\h]_{NR}=0.$ We define a linear map $\dM=:[\lambda\mu_\h,\cdot]_{NR}$ on the graded space $C^*(\g\oplus\h,\g\oplus\h)$. For all $g\in\Hom^n(S(\h),\g)$, we have
\begin{eqnarray*}
&&(\dM g)_p\Big((x_1,v_1),\cdots,(x_p,v_p)\Big)=[\lambda\mu_\h,g]_{NR}^p\Big((x_1,v_1),\cdots,(x_p,v_p)\Big)\\
&=&\Big(\lambda\mu_\h\circ g_{p-1}-(-1)^ng_{p-1}\circ \lambda\mu_\h\Big)\Big((x_1,v_1),\cdots,(x_p,v_p)\Big)\\
&=&\sum_{1\le i< j\le p}(-1)^{n+1+v_i(v_1+\cdots+v_{i-1})+v_j(v_1+\cdots+v_{j-1})+v_iv_j}
\Big(g_{p-1}(\lambda[v_i,v_j]_\h,v_1,\cdots,\hat{v}_i,\cdots,\hat{v}_j,\cdots,v_p),0\Big).
\end{eqnarray*}
Thus $\dM$ is closed on the subspace $sC^*(\h,\g)$, and is given by \eqref{diff}. By $[\lambda\mu_\h,\lambda\mu_\h]_{NR}=0$, we obtain that $\dM^2=0.$ Moreover, by $[\mu_\g+\rho,\lambda\mu_\h]_{NR}=0$, we deduce that $\dM$ is a derivation of the \gla~$(sC^*(\h,\g),\Courant{\cdot,\cdot})$. Therefore,  $(sC^*(\h,\g),\Courant{\cdot,\cdot},\dM)$ is a \dgla.
\end{proof}

Homotopy  $\huaO$-operators of weight $\lambda$ can be characterized as Maurer-Cartan elements of the above \dgla. Note that an element $T=\sum_{i=0}^{+\infty}T_i\in \Hom(S(\h),\g)$ is of degree $0$ if and only if the corresponding element $T\in s\Hom(S(\h),\g)$ is of degree $1$.

\begin{thm}\label{hmotopy-o-operator-dgla}
Let $\rho$ be an action of a \sgla~ $(\g,[\cdot,\cdot]_\g)$ on a \sgla~ $(\h,[\cdot,\cdot]_\h)$. A degree $0$ element $T=\sum_{i=0}^{+\infty}T_i\in \Hom(S(\h),\g)$ is a homotopy $\huaO$-operator of weight $\lambda$ on $\g$ with respect to the action $\rho$ if and only if $T=\sum_{i=0}^{+\infty}T_i$ is a Maurer-Cartan element of the \dgla $(sC^*(\h,\g),\Courant{\cdot,\cdot},\dM)$, i.e.
$$\dM T+\half\Courant{T,T}=0.$$
\end{thm}
\begin{proof}
For a degree $0$ element $T=\sum_{i=0}^{+\infty}T_i$ of the graded vector space $C^*(\h,\g)$, we write
$$
\dM T+\half\Courant{T,T}=\sum_i(\dM T+\half\Courant{T,T})_i\quad \mbox{with}\quad (\dM T+\half\Courant{T,T})_i\in\Hom(S^i(\h),\g).
$$
By Theorem \ref{dgla-deforamtion-homotopy}, we have
\begin{eqnarray*}
&&(\dM T+\half\Courant{T,T})_p(v_1,\cdots,v_p)\\
&=&\sum_{1\le i< j\le p}(-1)^{1+v_i(v_1+\cdots+v_{i-1})+v_j(v_1+\cdots+v_{j-1})+v_iv_j}T_{p-1}(\lambda[v_i,v_j]_\h,v_1,\cdots,\hat{v}_i,\cdots,\hat{v}_j,\cdots,v_p)\\
&&-\half\sum_{k+l=p+1}\sum_{\sigma\in \mathbb S_{(l,1,p-l-1)}}\varepsilon(\sigma)T_{k-1}\Big(\rho\big(T_l(v_{\sigma(1)},\cdots,v_{\sigma(l)})\big)v_{\sigma(l+1)},v_{\sigma(l+2)},\cdots,v_{\sigma(p)}\Big)\\
&&-\half\sum_{k+l=p+1}\sum_{\sigma\in \mathbb S_{(k-1,1,p-k)}}\varepsilon(\sigma)T_l\Big(\rho\big(T_{k-1}(v_{\sigma(1)},\cdots,v_{\sigma(k-1)})\big)v_{\sigma(k)}, v_{\sigma(k+1)},\cdots,v_{\sigma(p)}\Big)\\
&&+\half\sum_{k+l=p+1}\sum_{\sigma\in \mathbb S_{(k-1,l)}}\varepsilon(\sigma)[T_{k-1}(v_{\sigma(1)},\cdots,v_{\sigma(k-1)}),T_l(v_{\sigma(k)},\cdots,v_{\sigma(p)})]_\g\\
&=&-\sum_{1\le i< j\le p}(-1)^{v_i(v_1+\cdots+v_{i-1})+v_j(v_1+\cdots+v_{j-1})+v_iv_j}T_{p-1}(\lambda[v_i,v_j]_\h,v_1,\cdots,\hat{v}_i,\cdots,\hat{v}_j,\cdots,v_p)\\
&&-\sum_{k+l=p+1}\sum_{\sigma\in \mathbb S_{(l,1,p-l-1)}}\varepsilon(\sigma)T_{k-1}\Big(\rho\big(T_l(v_{\sigma(1)},\cdots,v_{\sigma(l)})\big)v_{\sigma(l+1)},v_{\sigma(l+2)},\cdots,v_{\sigma(p)}\Big)\\
&&+\half\sum_{k+l=p+1}\sum_{\sigma\in \mathbb S_{(k-1,l)}}\varepsilon(\sigma)[T_{k-1}(v_{\sigma(1)},\cdots,v_{\sigma(k-1)}),T_l(v_{\sigma(k)},\cdots,v_{\sigma(p)})]_\g.
\end{eqnarray*}
Thus,  $T=\sum_{i=0}^{+\infty}T_i\in \Hom(S(\h),\g)$ is a  homotopy $\huaO$-operator of weight $\lambda$ on $\g$ with respect to the action $\rho$ if and only if $T=\sum_{i=0}^{+\infty}T_i$ is a Maurer-Cartan element of the \dgla~ $(sC^*(\h,\g),\Courant{\cdot,\cdot},\dM)$.
\end{proof}

We note that a Lie algebra is a \sgla~ concentrated at degree $-1$. Moreover, a Lie algebra action is the same as an action of the \sgla~ on a graded vector space concentrated at degree $-1$. Therefore, we have the following corollary.

\begin{cor}
Let $\rho:\g\lon\Der(\h)$ be an action of a Lie algebra $\g$ on a Lie algebra $\h$. Then a linear map $T:\h\lon\g$  is an $\huaO$-operator of weight $\lambda$ on $\g$ with respect to the action $\rho$ if and only if $T$ is a Maurer-Cartan element of the \dgla~ $(\oplus_{n=0}^{\dim(\h)}\Hom(\wedge^{n}\h,\g),\Courant{\cdot,\cdot},\dM)$, where the differential $\dM:\Hom(\wedge^{n}\h,\g)\lon\Hom(\wedge^{n+1}\h,\g)$ is given by
\begin{eqnarray*}
(\dM g) ( v_1,\cdots, v_{n+1} )=\sum_{1\le i< j\le n+1}(-1)^{n+i+j-1}g(\lambda[v_i,v_j]_{\h},v_1,\cdots,\hat{v}_i,\cdots,\hat{v}_{j},\cdots,v_{n+1}),
\end{eqnarray*}
for all $g\in C^n(\h,\g)$ and $v_1,\cdots, v_{n+1} \in\h$, and the graded Lie bracket  $$\Courant{\cdot,\cdot}: \Hom(\wedge^n\h,\g)\times \Hom(\wedge^m\h,\g)\longrightarrow \Hom(\wedge^{m+n}\h,\g)$$ is given by
\begin{eqnarray*}
&&\Courant{g_1,g_2} ( v_1,\cdots, v_{m+n} )\\
&:=&-\sum_{\sigma\in \mathbb S_{(m,1,n-1)}}(-1)^{\sigma}g_1\Big(\rho\big(g_2(v_{\sigma(1)},\cdots,v_{\sigma(m)})\big)v_{\sigma(m+1)}, v_{\sigma(m+2)},\cdots,v_{\sigma(m+n)}\Big)\\
&&+(-1)^{mn}\sum_{\sigma\in \mathbb S_{(n,1,m-1)}}(-1)^{\sigma}g_2\Big(\rho\big(g_1(v_{\sigma(1)},\cdots,v_{\sigma(n)})\big)v_{\sigma(n+1)}, v_{\sigma(n+2)},\cdots,v_{\sigma(m+n)}\Big)\\
&&-(-1)^{mn}\sum_{\sigma\in \mathbb S_{(n,m)}}(-1)^{\sigma}[g_1(v_{\sigma(1)},\cdots,v_{\sigma(n)}),g_2(v_{\sigma(n+1)},\cdots,v_{\sigma(m+n)})]_{\g}
\end{eqnarray*}
for all  $g_1\in \Hom(\wedge^n\h,\g),~g_2\in \Hom(\wedge^m\h,\g)$ and $v_1,\cdots, v_{m+n} \in\h.$
\end{cor}

If the Lie algebra $\h$ is abelian in the above corollary, we recover the  \gla~ that controls the deformations of $\huaO$-operators of weight $0$ given in \cite[Proposition~2.3]{TBGS}.

\section{\Opt homotopy post-Lie algebras}
\label{sec:post}
In this section, we first recall the notion of a post-Lie algebra, and then give the definition of an \opt homotopy post-Lie algebra as a variation of a homotopy post-Lie algebra. We construct a \dgla~ and show that \opt homotopy post-Lie algebras can be characterized as its Maurer-Cartan elements to justify the notion. We also show that  \opt homotopy post-Lie  algebras naturally arise from homotopy $\huaO$-operators of weight 1.

\begin{defi} (\cite{Val})
A {\bf post-Lie algebra} $(\g,[\cdot,\cdot]_\g,\rhd)$ consists of a Lie algebra $(\g,[\cdot,\cdot]_\g)$ and a binary product $\rhd:\g\otimes\g\lon\g$ such that
\begin{eqnarray}
\label{Post-1}x\rhd[y,z]_\g&=&[x\rhd y,z]_\g+[y,x\rhd z]_\g,\\
\label{Post-2}[x,y]_\g\rhd z&=&a_\rhd(x,y,z)-a_\rhd(y,x,z),
\end{eqnarray}
here $a_{\rhd}(x,y,z):=x\rhd(y\rhd z)-(x\rhd y)\rhd z $ and $x,y,z\in \g.$
\end{defi}

Define $L_\rhd:\g\lon \gl(\g)$ by $L_\rhd(x)(y)=x\rhd y$. Then by \eqref{Post-1},    $L_\rhd$ is a linear map from $\g$ to $\Der(\g)$. In the sequel, we will say that $ L_\rhd$ is a post-Lie algebra structure on the Lie algebra  $(\g,[\cdot,\cdot]_\g)$.
\begin{rmk}
Let $(\g,[\cdot,\cdot]_\g,\rhd)$ be a post-Lie algebra. If the Lie bracket $[\cdot,\cdot]_\g=0$, then $(\g,\rhd)$ becomes a pre-Lie algebra. Thus,  a post-Lie algebra can be viewed as a nonabelian version of a pre-Lie algebra. See \cite{Burde16,Burde19} for the classifications of post-Lie algebras on certain Lie algebras, and \cite{Munthe-Kaas-Lundervold} for applications of post-Lie algebras in numerical integration.
\end{rmk}

The following well-known result is a special case of splitting of operads~\cite{BBGN,PBG,Val}.
\begin{pro}\label{double-Lie}
Let $(\g,[\cdot,\cdot]_\g,\rhd)$ be a post-Lie algebra. Then    the bracket $[\cdot,\cdot]_C$ defined by
\begin{eqnarray}
[x,y]_C:=x\rhd y-y\rhd x+[x,y]_\g,\quad\forall x,y\in\g,
\end{eqnarray}
is a Lie bracket.
\end{pro}

We denote this Lie algebra by $\g^C$ and call it  the {\bf sub-adjacent Lie algebra} of $(\g,[\cdot,\cdot]_\g,\rhd)$.

\begin{defi}
An {\bf \opt homotopy post-Lie algebra} is a \sgla~ $(\g,[\cdot,\cdot]_\g)$ equipped with a collection $(k\ge 1)$ of linear maps $\oprn_k:\otimes^k \g\lon \g$ of degree $1$ satisfying, for every collection of
homogeneous elements $x_1,\cdots,x_n,x_{n+1}\in \g$,
\begin{itemize}
\item[\rm(i)]
({\bf graded symmetry}) for every $\sigma\in\mathbb S_{n-1},~n\ge1$,
\begin{eqnarray}
\label{homotopy-post-1}\oprn_n(x_{\sigma(1)},\cdots,x_{\sigma(n-1)},x_n)=\varepsilon(\sigma)\oprn_n(x_1,\cdots,x_{n-1},x_n),
\end{eqnarray}
\item[\rm(ii)]
({\bf graded derivation}) for all $n\ge 1$,
\begin{eqnarray}
\label{homotopy-post-2}\nonumber\oprn_n(x_1,\cdots,x_{n-1},[x_n,x_{n+1}]_\g)
&=&(-1)^{x_1+\cdots+x_{n-1}+1}[\oprn_n(x_1,\cdots,x_{n-1},x_n),x_{n+1}]_\g\\
&&+(-1)^{(x_1+\cdots+x_{n-1}+1)(x_n+1)}[x_n,\oprn_n(x_1,\cdots,x_{n-1},x_{n+1})]_\g,
\end{eqnarray}
\item[\rm(iii)] for all $n\ge 1$,
\begin{eqnarray}
\nonumber &&\sum_{i+j=n+1\atop i\ge1,j\geq2}\sum_{\sigma\in\mathbb S_{(i-1,1,j-2)}}\varepsilon(\sigma)\oprn_j(\oprn_i(v_{\sigma(1)},\cdots,v_{\sigma(i-1)},v_{\sigma(i)}), v_{\sigma(i+1)},\cdots,v_{\sigma(n-1)},v_{n})\\
\label{homotopy-post-3}
&&+\sum_{i+j=n+1\atop i\geq1,j\geq1}\sum_{\sigma\in\mathbb S_{(j-1,i-1)}}(-1)^{\alpha}\varepsilon(\sigma)\oprn_j(v_{\sigma(1)},\cdots,v_{\sigma(j-1)},\oprn_i( v_{\sigma(j)},\cdots,v_{\sigma(n-1)},v_{n})) \\
&=&\nonumber\sum_{1\le i< j\le n-1}(-1)^{\beta}\oprn_{n-1}([x_i,x_j]_\g,x_1,\cdots,\hat{x}_i,\cdots,\hat{x}_j,\cdots,x_n),
\end{eqnarray}
where $\beta=x_i(x_1+\cdots+x_{i-1})+x_j(x_1+\cdots+x_{j-1})+x_ix_j+1$ and $\alpha=x_{\sigma(1)}+x_{\sigma(2)}+\cdots+x_{\sigma(j-1)}$.
\end{itemize}
\end{defi}

The notion of a pre-Lie$_\infty$ algebra was introduced in \cite{CL}. See \cite{Mer} for more applications of pre-Lie$_\infty$ algebras in geometry. Recall that a {\bf pre-Lie$_\infty$ algebra} is a graded vector space $V$ equipped with a collection of linear maps $\oprn_k:\otimes^k V\lon V, k\ge 1,$  of degree $1$ with the property that, for any homogeneous elements $v_1,\cdots,v_n\in V$, we have
\begin{itemize}\item[\rm(i)]
{\bf (graded symmetry)} for every $\sigma\in\mathbb S_{n-1}$,
\begin{eqnarray*}
\oprn_n(v_{\sigma(1)},\cdots,v_{\sigma(n-1)},v_n)=\varepsilon(\sigma)\oprn_n(v_1,\cdots,v_{n-1},v_n),
\end{eqnarray*}
\item[\rm(ii)] for all $n\ge 1$,
\begin{eqnarray*}
&&\sum_{i+j=n+1\atop i\ge1,j\geq2}\sum_{\sigma\in\mathbb S_{(i-1,1,j-2)}}\varepsilon(\sigma)\oprn_j(\oprn_i(v_{\sigma(1)},\cdots,v_{\sigma(i-1)},v_{\sigma(i)}), v_{\sigma(i+1)},\cdots,v_{\sigma(n-1)},v_{n})\\
&&+\sum_{i+j=n+1\atop i\geq1,j\geq1}\sum_{\sigma\in\mathbb S_{(j-1,i-1)}}(-1)^{\alpha}\varepsilon(\sigma)\oprn_j(v_{\sigma(1)},\cdots,v_{\sigma(j-1)},\oprn_i( v_{\sigma(j)},\cdots,v_{\sigma(n-1)},v_{n}))=0,
\end{eqnarray*}
where $\alpha=v_{\sigma(1)}+v_{\sigma(2)}+\cdots+v_{\sigma(j-1)}$.
\end{itemize}
It is straightforward to see that an \opt homotopy post-Lie algebra $(\g,[\cdot,\cdot]_\g,\{\oprn_k\}_{k=1}^\infty)$ reduces to a pre-Lie$_\infty$ algebra when $(\g,[\cdot,\cdot]_\g)$ is an abelian \sgla.

\begin{rmk}
The explicit definition of a homotopy post-Lie algebra is still not known though its operad should be the minimal model of the post-Lie operad by general construction as noted in the introduction~\cite{LV}. Since  a Rota-Baxter Lie algebra of weight $1$ gives  a post-Lie algebra, we expect that a Rota-Baxter homotopy Lie algebra of weight $1$ (a homotopy Lie algebra with a Rota-Baxter operator of weight 1) induces a homotopy post-Lie algebra. Based on preliminary computations in this regards, there should be more terms in the right hand side of \eqref{homotopy-post-3} in the definition of a homotopy post-Lie algebra, suggesting that the homotopy post-Lie algebra is different from the \opt homotopy post-Lie algebra just defined. So we have chosen a different name for distinction.
\label{rk:postlie}
\end{rmk}

Now we  construct the  \dgla ~ that characterizes \opt homotopy post-Lie algebras as Maurer-Cartan elements. Let $(\h,[\cdot,\cdot]_\h)$ be a \sgla. Denote by
$$
\bar{\huaC}^n(\h,\h)=\Hom^n(S(\h),s^{-1}\Der(\h)),\quad \bar{\huaC}^*(\h,\h):=\oplus_{n\in\mathbb Z} \bar{\huaC}^n(\h,\h).
$$

\begin{itemize}
\item[$\bullet$]Define a graded linear map $\partial:\bar{\huaC}^n(\h,\h)\lon\bar{\huaC}^{n+1}(\h,\h)$ by
\begin{eqnarray*}
(\partial \beta)_p (v_1,\cdots,v_{p})=\sum_{1\le i< j\le s}(-1)^{n+1+v_i(v_1+\cdots+v_{i-1})+v_j(v_1+\cdots+v_{j-1})+v_iv_j}\beta_{p-1}([v_i,v_j]_{\h},v_1,\cdots,\hat{v}_i,\cdots,\hat{v}_{j},\cdots,v_{p}),
\end{eqnarray*}

\item[$\bullet$]Define a graded bracket operation
 $$[\cdot,\cdot]^c:\bar{\huaC}^m(\h,\h)\times \bar{\huaC}^n(\h,\h)\longrightarrow \bar{\huaC}^{m+n+1}(\h,\h)$$
by
\begin{eqnarray*}
&&[\alpha,\beta]^c_p ( v_1,\cdots, v_{p} )\\
&=&-\sum_{k+l=p+1}\sum_{\sigma\in \mathbb S_{(l,1,p-l-1)}}\varepsilon(\sigma)\alpha_{k-1}\Big(s\big(\beta_l(v_{\sigma(1)},\cdots,v_{\sigma(l)})\big)v_{\sigma(l+1)},v_{\sigma(l+2)},\cdots,v_{\sigma(p)}\Big)\\
&&+(-1)^{(m+1)(n+1)}\sum_{k+l=p+1}\sum_{\sigma\in \mathbb S_{(k-1,1,p-k)}}\varepsilon(\sigma)\beta_l\Big(s\big(\alpha_{k-1}(v_{\sigma(1)},\cdots,v_{\sigma(k-1)})\big)v_{\sigma(k)}, v_{\sigma(k+1)},\cdots,v_{\sigma(p)}\Big)\\
&&-\sum_{k+l=p+1}\sum_{\sigma\in \mathbb S_{(k-1,l)}}(-1)^{(n+1)(v_{\sigma(1)}+\cdots+v_{\sigma(k-1)})}\varepsilon(\sigma)s^{-1}[s\alpha_{k-1}(v_{\sigma(1)},\cdots,v_{\sigma(k-1)}),s\beta_l(v_{\sigma(k)},\cdots,v_{\sigma(p)})]
%&&-\sum_{k+l=s+1}\sum_{\sigma\in \mathbb S_{(s-l,l)}}(-1)^{(n+1)(v_{\sigma(1)}+\cdots+v_{\sigma(s-l)})}\varepsilon(\sigma)s^{-1}[s\alpha_{k-1}(v_{\sigma(1)},\cdots,v_{\sigma(s-l)}),s\beta_l(v_{\sigma(s-l+1)},\cdots,v_{\sigma(s)})]
\end{eqnarray*}
for all $\alpha=\sum_i\alpha_i\in\bar{\huaC}^m(\h,\h),$ $\beta=\sum_i\beta_i\in \bar{\huaC}^n(\h,\h)$ with $\alpha_i,~\beta_i\in\Hom(S^i(\h),s^{-1}\Der(\h))$ and $v_1,\cdots, v_{p} \in\h.$ Here we write $\partial\beta=\sum_i(\partial \beta)_i$ with $(\partial \beta)_i\in\Hom(S^i(\h),s^{-1}\Der(\h))$ and $[\alpha,\beta]^c=\sum_i[\alpha,\beta]^c_i$ with $[\alpha,\beta]^c_i\in\Hom(S^i(\h),s^{-1}\Der(\h))$.
\end{itemize}

\begin{thm}\label{homotopy-post-dg-lie}
 With  above notations, $(s\bar{\huaC}^*(\h,\h),[\cdot,\cdot]^c,\partial)$ is a \dgla. Moreover, its Maurer-Cartan elements are precisely  the \opt homotopy post-Lie algebra structures on the \sgla~ $(\h,[\cdot,\cdot]_\h)$.
\end{thm}

\begin{proof}
Setting $\lambda=1,\g=s^{-1}\Der(\h)$ and $\rho=s$ in Theorem~\ref{dgla-deforamtion-homotopy}, we find that  $(s\bar{\huaC}^*(\h,\h),[\cdot,\cdot]^c,\partial)$ is a \dgla.

Let $L=\sum_{i=0}^{+\infty}L_i\in \Hom^0(S(\h),s^{-1}\Der(\h))$ with $L_i\in\Hom(S^i(\h),s^{-1}\Der(\h))$ be a Maurer-Cartan element of the \dgla~ $(s\bar{\huaC}^*(\h,\h),[\cdot,\cdot]^c,\partial)$. We define a collection of linear maps $\oprn_k:\otimes^k \h\lon \h$  $(k\ge 1)$ of degree $1$ by
\begin{eqnarray}
\label{homotopy-post-lie}\oprn_k(v_1,\cdots,v_k):=\big(sL_{k-1}(v_1,\cdots,v_{k-1})\big)v_k,\quad\forall v_1\cdots,v_k\in\h.
\end{eqnarray}
By $L_{n-1}\in\Hom^0(S^{n-1}(\h),s^{-1}\Der(\h))$, we obtain \eqref{homotopy-post-1} and \eqref{homotopy-post-2}. Moreover, write $\partial L+\half[L,L]^c=\sum_i(\partial L+\half[L,L]^c)_{i}$ with $(\partial L+\half[L,L]^c)_{i}\in\Hom(S^i(\h),s^{-1}\Der(\h))$. Then for all $v_1,\cdots,v_n\in\h$, by a similar computation as in the proof of Theorem~\ref{hmotopy-o-operator-dgla}, we have
\begin{eqnarray*}
&&\Big(s\big((\partial L+\half[L,L]^c)_{n-1}\big)(v_1,\cdots,v_{n-1})\Big)v_n\\
&=&-\sum_{1\le i< j\le n-1}(-1)^{v_i(v_1+\cdots+v_{i-1})+v_j(v_1+\cdots+v_{j-1})+v_iv_j}sL_{n-2}([v_i,v_j]_{\h},v_1,\cdots,\hat{v}_i,\cdots,\hat{v}_{j},\cdots,v_{n-1})v_n\\
&&-\sum_{k+l=n}\sum_{\sigma\in \mathbb S_{(l,1,n-l-2)}}\varepsilon(\sigma)sL_{k-1}\Big(s\big(L_l(v_{\sigma(1)},\cdots,v_{\sigma(l)})\big)v_{\sigma(l+1)},v_{\sigma(l+2)},\cdots,v_{\sigma(n-1)}\Big)v_n\\
&&-\half\sum_{k+l=n}\sum_{\sigma\in \mathbb S_{(k-1,l)}}(-1)^{v_{\sigma(1)}+\cdots+v_{\sigma(k-1)}}\varepsilon(\sigma)[sL_{k-1}(v_{\sigma(1)},\cdots,v_{\sigma(k-1)}),sL_l(v_{\sigma(k)},\cdots,v_{\sigma(n-1)})]v_n\\
&=&-\sum_{1\le i< j\le n-1}(-1)^{v_i(v_1+\cdots+v_{i-1})+v_j(v_1+\cdots+v_{j-1})+v_iv_j}\oprn_{n-1}([v_i,v_j]_{\h},v_1,\cdots,\hat{v}_i,\cdots,\hat{v}_{j},\cdots,v_{n-1},v_n)\\
&&-\sum_{k+l=n}\sum_{\sigma\in \mathbb S_{(l,1,n-l-2)}}\varepsilon(\sigma)\oprn_{k}\big(\oprn_{l+1}(v_{\sigma(1)},\cdots,v_{\sigma(l)},v_{\sigma(l+1)}),v_{\sigma(l+2)},\cdots,v_{\sigma(n-1)},v_n\big)\\
&&-\half\sum_{k+l=n}\sum_{\sigma\in \mathbb S_{(k-1,l)}}(-1)^{v_{\sigma(1)}+\cdots+v_{\sigma(k-1)}}\varepsilon(\sigma)\oprn_{k}\big(v_{\sigma(1)},\cdots,v_{\sigma(k-1)},\oprn_{l+1}(v_{\sigma(k)},\cdots,v_{\sigma(n-1)},v_n)\big)\\
&&-\half\sum_{k+l=n}\sum_{\sigma\in \mathbb S_{(k-1,l)}}(-1)^{(v_{\sigma(1)}+\cdots+v_{\sigma(k-1)}+1)(v_{\sigma(k)}+\cdots+v_{\sigma(n-1)})}\varepsilon(\sigma)\oprn_{l+1}\big(v_{\sigma(k)},\cdots,v_{\sigma(n-1)},\oprn_{k}(v_{\sigma(1)},\cdots,v_{\sigma(k-1)},v_n)\big).
\end{eqnarray*}
For any $\sigma\in\mathbb S_{(k-1,l)}$, we define $\tau\in\mathbb S_{(l,k-1)}$ by
\[
\tau(i):=\left\{
\begin{array}{ll}
\sigma(i+k-1), & 1\le i\le l;\\
\sigma(i-l), & l+1\le i\le n-1.
\end{array}
\right.
\]
Thus we have $\varepsilon(\tau;v_1,\cdots,v_{n-1})=\varepsilon(\sigma;v_1,\cdots,v_{n-1})(-1)^{(v_{\sigma(1)}+\cdots+v_{\sigma(k-1)})(v_{\sigma(k)}+\cdots+v_{\sigma(n-1)})}$.
Applying the bijection between $\mathbb S_{(k-1,l)}$ and $\mathbb S_{(l,k-1)}$, we obtain
\begin{eqnarray}
\nonumber&&-\half\sum_{k+l=n}\sum_{\sigma\in \mathbb S_{(k-1,l)}}(-1)^{(v_{\sigma(1)}+\cdots+v_{\sigma(k-1)}+1)(v_{\sigma(k)}+\cdots+v_{\sigma(n-1)})}\varepsilon(\sigma)\oprn_{l+1}\big(v_{\sigma(k)},\cdots,v_{\sigma(n-1)},\oprn_{k}(v_{\sigma(1)},\cdots,v_{\sigma(k-1)},v_n)\big)\\
\nonumber&&=-\half\sum_{k+l=n}\sum_{\tau\in \mathbb S_{(l,k-1)}}(-1)^{v_{\tau(1)}+\cdots+v_{\tau(l)}}\varepsilon(\tau)\oprn_{l+1}\big(v_{\tau(1)},\cdots,v_{\tau(l)},\oprn_{k}(v_{\tau(l+1)},\cdots,v_{\tau(n-1)},v_n)\big).
\end{eqnarray}
Therefore, we have
\begin{eqnarray*}
&&\Big(s\big((\partial L+\half[L,L]^c)_{n-1}\big)(v_1,\cdots,v_{n-1})\Big)v_n\\
&=&-\sum_{1\le i< j\le n-1}(-1)^{v_i(v_1+\cdots+v_{i-1})+v_j(v_1+\cdots+v_{j-1})+v_iv_j}\oprn_{n-1}([v_i,v_j]_{\h},v_1,\cdots,\hat{v}_i,\cdots,\hat{v}_{j},\cdots,v_{n-1},v_n)\\
&&-\sum_{k+l=n}\sum_{\sigma\in \mathbb S_{(l,1,n-l-2)}}\varepsilon(\sigma)\oprn_{k}\big(\oprn_{l+1}(v_{\sigma(1)},\cdots,v_{\sigma(l)},v_{\sigma(l+1)}),v_{\sigma(l+2)},\cdots,v_{\sigma(n-1)},v_n\big)\\
&&-\sum_{k+l=n}\sum_{\sigma\in \mathbb S_{(k-1,l)}}(-1)^{v_{\sigma(1)}+\cdots+v_{\sigma(k-1)}}\varepsilon(\sigma)\oprn_{k}\big(v_{\sigma(1)},\cdots,v_{\sigma(k-1)},\oprn_{l+1}(v_{\sigma(k)},\cdots,v_{\sigma(n-1)},v_n)\big),
\end{eqnarray*}
which implies \eqref{homotopy-post-3}. Thus,   $\{\oprn_k\}_{k=1}^\infty$  is an \opt homotopy post-Lie algebra structure on the \sgla~ $(\h,[\cdot,\cdot]_\h)$.
\end{proof}

When the \sgla~ $\h$ reduces to a usual Lie algebra, we characterize post-Lie algebra structures on $\h$ as Maurer-Cartan elements. See~\cite{Burde16,Burde19,GLBJ,PLBG} for classifications of post-Lie algebras on some specific Lie algebras.

\begin{cor}\label{cor:dgla-post-lie}
 Let $(\h,[\cdot,\cdot]_\h)$ be a Lie algebra. Denote by
$
\bar{\huaC}^n(\h,\h):=\Hom(\wedge^n\h,\Der(\h))$ and $ \bar{\huaC}^*(\h,\h):=\oplus_{n\geq 0} \bar{\huaC}^n(\h,\h).
$ Then $(\bar{\huaC}^*(\h,\h),[\cdot,\cdot]^c,\partial)$ is a \dgla, where  the differential $\partial:\bar{\huaC}^n(\h,\h)\lon\bar{\huaC}^{n+1}(\h,\h)$ is given by
\begin{eqnarray*}
(\partial \alpha) ( u_1,\cdots, u_{n+1} )=\sum_{1\le i< j\le n+1}(-1)^{n+i+j-1}\alpha([u_i,u_j]_{\h},u_1,\cdots,\hat{u}_i,\cdots,\hat{u}_{j},\cdots,u_{n+1}),\,\,~
\end{eqnarray*}
 and the graded Lie bracket
$
[\cdot,\cdot]^c: \bar{\huaC}^n(\h,\h)\times \bar{\huaC}^m(\h,\h)\lon\bar{\huaC}^{m+n}(\h,\h)
$ is given by
\begin{eqnarray*}
&&[\alpha,\beta]^c(u_1,\cdots,u_{m+n})\\
&=&-\sum_{\sigma\in\mathbb S_{(m,1,n-1)}}(-1)^{\sigma}\alpha(\beta(u_{\sigma(1)},\cdots,u_{\sigma(m)})u_{\sigma(m+1)},u_{\sigma(m+2)},\cdots,u_{\sigma(m+n)})\\
&&+(-1)^{mn}\sum_{\sigma\in\mathbb S_{(n,1,m-1)}}(-1)^{\sigma}\beta(\alpha(u_{\sigma(1)},\cdots,u_{\sigma(n)})u_{\sigma(n+1)},u_{\sigma(n+2)},\cdots,u_{\sigma(m+n)})\\
&&-(-1)^{mn}\sum_{\sigma\in\mathbb S_{(n,m)}}(-1)^{\sigma}[\alpha(u_{\sigma(1)},\cdots,u_{\sigma(n)}),\beta(u_{\sigma(n+1)},\cdots,u_{\sigma(m+n)})],
\end{eqnarray*}
for all $\alpha\in \Hom(\wedge^{n}\h,\Der(\h))$, $\beta\in \Hom(\wedge^{m}\h,\Der(\h))$ and $u_1,\cdots, u_{m+n}\in\h$.

Moreover, $L_\rhd:\h\lon\Der(\h)$ defines a post-Lie algebra structure on the Lie algebra $(\h,[\cdot,\cdot]_\h)$ if and only if $L_\rhd$ is a Maurer-Cartan element of the \dgla~ $(\bar{\huaC}^*(\h,\h),[\cdot,\cdot]^c,\partial)$.
\end{cor}

When the \sgla~ $\frkh$ is abelian, we characterize pre-Lie$_\infty$-algebras as Maurer-Cartan elements, which was originally given in \cite{CL}.

\begin{cor}
Let $V$ be a $\mathbb Z$-graded vector space. Denote
$$
\bar{\huaC}^n(V,V):=\Hom^n(S(V),\gl(V)),\quad \bar{\huaC}^*(V,V):=\oplus_{n\in\mathbb Z} \bar{\huaC}^n(V,V).
$$
Then $(\bar{\huaC}^*(V,V),[\cdot,\cdot]^c)$ is a \gla, where  the graded Lie bracket
 $[\cdot,\cdot]^c:\bar{\huaC}^m(V,V)\times \bar{\huaC}^n(V,V)\longrightarrow \bar{\huaC}^{m+n}(V,V)$ is given
by
\begin{eqnarray*}
&&[\alpha,\beta]^c_p ( v_1,\cdots, v_{p} )\\
&=&-\sum_{k+l=p+1}\sum_{\sigma\in \mathbb S_{(l,1,p-l-1)}}\varepsilon(\sigma)\alpha_{k-1}\Big(\beta_l(v_{\sigma(1)},\cdots,v_{\sigma(l)})v_{\sigma(l+1)},v_{\sigma(l+2)},\cdots,v_{\sigma(p)}\Big)\\
&&+(-1)^{mn}\sum_{k+l=p+1}\sum_{\sigma\in \mathbb S_{(k-1,1,p-k)}}\varepsilon(\sigma)\beta_l\Big(\alpha_{k-1}(v_{\sigma(1)},\cdots,v_{\sigma(k-1)})v_{\sigma(k)}, v_{\sigma(k+1)},\cdots,v_{\sigma(p)}\Big)\\
&&-\sum_{k+l=p+1}\sum_{\sigma\in \mathbb S_{(k-1,l)}}(-1)^{n(v_{\sigma(1)}+\cdots+v_{\sigma(k-1)})}\varepsilon(\sigma)[\alpha_{k-1}(v_{\sigma(1)},\cdots,v_{\sigma(k-1)}),\beta_l(v_{\sigma(k)},\cdots,v_{\sigma(p)})],
\end{eqnarray*}
for all $\alpha=\sum_i\alpha_i\in\bar{\huaC}^m(V,V),~\beta=\sum\beta_i\in\bar{\huaC}^n(V,V)$ with $\alpha_i, \beta_i\in\Hom(S^i(V),V)$.

Moreover, $L=\sum_{i=0}^{\infty}L_i\in\Hom^1(S(\h),\gl(\h))$ defines a pre-Lie$_\infty$ algebra structure by
\begin{eqnarray}
\label{homotopy-pre-lie}\oprn_k(v_1,\cdots,v_k):=L_{k-1}(v_1,\cdots,v_{k-1})v_k,\quad\forall v_1\cdots,v_k\in\h,
\end{eqnarray}
on the graded vector space $V$ if and only if  $L=\sum_{i=0}^{\infty}L_i$ is a  Maurer-Cartan element of the \gla~ $(\bar{\huaC}^*(V,V),[\cdot,\cdot]^c)$.
\end{cor}

In the above corollary, if the graded vector space $V$ reduces to a usual vector space, we characterize pre-Lie algebra structures as Maurer-Cartan elements. See \cite{CL,Nijenhuis,WBLS} for more details.

It is known that post-Lie algebras naturally arise from $\huaO$-operators of weight $1$ as follows.

\begin{pro}{\rm (\cite{BGN})}\label{rota-baxter-to-post-Lie}
Let $T:\h\to \g$ be an $\huaO$-operator of weight $1$. Then  $(\h,[\cdot,\cdot]_\h,\rhd)$ is a post-Lie algebra, where $\rhd$ is given by
\begin{eqnarray*}
u\rhd v=\rho(Tu)v,\,\,\,\,\forall u,v\in\h.
\end{eqnarray*}
\end{pro}

In the sequel, we generalize the above relation to homotopy $\huaO$-operators of weight $1$ and \opt homotopy post-Lie algebras.

Define a graded linear map $\Psi:C^*(\h,\g)\lon\bar{\huaC}^*(\h,\h)$ of degree $0$ by
$$
\Psi(f)=s^{-1}\circ\rho\circ f,\quad \forall f\in\Hom^m(S(\h),\g).
$$
Therefore, we have $\Psi(f)_k=s^{-1}\circ\rho\circ f_k.$
In the following, we set $\lambda=1$ in Theorem \ref{dgla-deforamtion-homotopy} for notational simplicity.

\begin{thm}\label{homo-dg-lie}
Let $(\g,[\cdot,\cdot]_\g)$ and $(\h,[\cdot,\cdot]_\h)$ be \sgla's and  $\rho:\g\lon\Der(\h)$ a \sgla~ action. Then $\Psi$ is a
homomorphism of \dgla's from
$(sC^*(\h,\g),\Courant{\cdot,\cdot},\dM)$  to
$(s\bar{\huaC}^*(\h,\h),[\cdot,\cdot]^c,\partial)$.
\end{thm}

\begin{proof}
 For all $f=\sum_i f_i\in\Hom^m(S(\h),\g),~ g=\sum_ig_i\in\Hom^n(S(\h),\g)$ with $f_i,~g_i\in\Hom(S^i(\h),\g)$, we write  $[\Psi(f),\Psi(g)]^c-\Psi(\Courant{f,g}) =\sum_i\Big([\Psi(f),\Psi(g)]^c-\Psi(\Courant{f,g})\Big)_i$. Since $s^{-1}\circ\rho:\g\lon s^{-1}\Der(\h)$ is a \sgla~ homomorphism,
we have
\begin{eqnarray*}
&&\Big([\Psi(f),\Psi(g)]^c-\Psi(\Courant{f,g})\Big)_p(v_1,\cdots,v_{p})\\
&=&-\sum_{k+l=p+1}\sum_{\sigma\in \mathbb S_{(l,1,p-l-1)}}\varepsilon(\sigma)\Psi(f)_{k-1}\Big(s\big(\Psi(g)_l(v_{\sigma(1)},\cdots,v_{\sigma(l)})\big)v_{\sigma(l+1)},
v_{\sigma(l+2)},\cdots,v_{\sigma(p)}\Big)\\
&&+(-1)^{(m+1)(n+1)}\sum_{k+l=p+1}\sum_{\sigma\in \mathbb S_{(k-1,1,p-k)}}\varepsilon(\sigma)\Psi(g)_l\Big(s\big(\Psi(f)_{k-1}(v_{\sigma(1)},\cdots,v_{\sigma(k-1)})\big)v_{\sigma(k)}, v_{\sigma(k+1)},\cdots,v_{\sigma(p)}\Big)\\
&&-\sum_{k+l=p+1}\sum_{\sigma\in \mathbb S_{(p-l,l)}}(-1)^{(n+1)(v_{\sigma(1)}+\cdots+v_{\sigma(p-l)})}\varepsilon(\sigma)s^{-1}[s\Psi(f)_{k-1}(v_{\sigma(1)},\cdots,v_{\sigma(p-l)}),s\Psi(g)_l(v_{\sigma(p-l+1)},\cdots,v_{\sigma(p)})]\\
&&+\sum_{k+l=p+1}\sum_{\sigma\in \mathbb S_{(l,1,p-l-1)}}\varepsilon(\sigma)s^{-1}\rho\Big(f_{k-1}\Big(\rho\big(g_l(v_{\sigma(1)},\cdots,v_{\sigma(l)})\big)v_{\sigma(l+1)},
v_{\sigma(l+2)},\cdots,v_{\sigma(p)}\Big)\Big)\\
\nonumber&&-(-1)^{(m+1)(n+1)}\sum_{k+l=p+1}\sum_{\sigma\in \mathbb S_{(k-1,1,p-k)}}\varepsilon(\sigma)s^{-1}\rho\Big(g_l\Big(\rho\big(f_{k-1}(v_{\sigma(1)},\cdots,v_{\sigma(k-1)})\big)v_{\sigma(k)}, v_{\sigma(k+1)},\cdots,v_{\sigma(p)}\Big)\Big)\\
\nonumber&&+\sum_{k+l=p+1}\sum_{\sigma\in \mathbb S_{(p-l,l)}}(-1)^{n(v_{\sigma(1)}+\cdots+v_{\sigma(p-l)})+m+1}\varepsilon(\sigma)s^{-1}\rho[f_{k-1}(v_{\sigma(1)},\cdots,v_{\sigma(p-l)}),
g_l(v_{\sigma(p-l+1)},\cdots,v_{\sigma(p)})]_\g\\
&=&-\sum_{k+l=p+1}\sum_{\sigma\in \mathbb S_{(l,1,p-l-1)}}\varepsilon(\sigma)s^{-1}\rho\Big(f_{k-1}\Big(\rho\big(g_l(v_{\sigma(1)},\cdots,v_{\sigma(l)})\big)v_{\sigma(l+1)},
v_{\sigma(l+2)},\cdots,v_{\sigma(p)}\Big)\Big)\\
&&+(-1)^{(m+1)(n+1)}\sum_{k+l=p+1}\sum_{\sigma\in \mathbb S_{(k-1,1,p-k)}}\varepsilon(\sigma)s^{-1}\rho\Big(g_l\Big(\rho(f_{k-1}(v_{\sigma(1)},\cdots,v_{\sigma(k-1)})\big)v_{\sigma(k)}, v_{\sigma(k+1)},\cdots,v_{\sigma(p)}\Big)\Big)\\
&&-\sum_{k+l=p+1}\sum_{\sigma\in \mathbb S_{(p-l,l)}}(-1)^{(n+1)(v_{\sigma(1)}+\cdots+v_{\sigma(p-l)})}\varepsilon(\sigma)s^{-1}[\rho\big(f_{k-1}(v_{\sigma(1)},\cdots,v_{\sigma(p-l)})\big),
\rho\big(g_l(v_{\sigma(p-l+1)},\cdots,v_{\sigma(p)})\big)]\\
&&+\sum_{k+l=p+1}\sum_{\sigma\in \mathbb S_{(l,1,p-l-1)}}\varepsilon(\sigma)s^{-1}\rho\Big(f_{k-1}\Big(\rho\big(g_l(v_{\sigma(1)},\cdots,v_{\sigma(l)})\big)v_{\sigma(l+1)},v_{\sigma(l+2)},
\cdots,v_{\sigma(p)}\Big)\Big)\\
\nonumber&&-(-1)^{(m+1)(n+1)}\sum_{k+l=p+1}\sum_{\sigma\in \mathbb S_{(k-1,1,p-k)}}\varepsilon(\sigma)s^{-1}\rho\Big(g_l\Big(\rho\big(f_{k-1}(v_{\sigma(1)},\cdots,v_{\sigma(k-1)})\big)v_{\sigma(k)}, v_{\sigma(k+1)},\cdots,v_{\sigma(p)}\Big)\Big)\\
\nonumber&&+\sum_{k+l=p+1}\sum_{\sigma\in \mathbb S_{(p-l,l)}}(-1)^{n(v_{\sigma(1)}+\cdots+v_{\sigma(p-l)})+m+1}\varepsilon(\sigma)s^{-1}\rho([f_{k-1}(v_{\sigma(1)},\cdots,v_{\sigma(p-l)}),
g_l(v_{\sigma(p-l+1)},\cdots,v_{\sigma(p)})]_\g)\\
&=&0.
\end{eqnarray*}
Moreover, for all $g\in\Hom^n(S(\h),\g)$, we have
\begin{eqnarray*}
&&\Big((\partial\circ\Psi-\Psi\circ\dM)(g)\Big)_p(v_1,\cdots,v_{p})\\
&=&\big(\partial(\Psi(g))\big)_p(v_1,\cdots,v_{p})-\big(\Psi(\dM g)\big)_p(v_1,\cdots,v_{p})\\
&=&\sum_{1\le i< j\le p}(-1)^{n+1+v_i(v_1+\cdots+v_{i-1})+v_j(v_1+\cdots+v_{j-1})+v_iv_j}(\Psi(g))_{p-1}([v_i,v_j]_{\h},v_1,\cdots,\hat{v}_i,\cdots,\hat{v}_{j},\cdots,v_{p})\\
&&-\sum_{1\le i< j\le p}(-1)^{n+1+v_i(v_1+\cdots+v_{i-1})+v_j(v_1+\cdots+v_{j-1})+v_iv_j}s^{-1}\rho\big(g_{p-1}([v_i,v_j]_{\h},v_1,\cdots,\hat{v}_i,\cdots,\hat{v}_{j},\cdots,
v_{p})\big)\\
&=&0.
\end{eqnarray*}
Thus, we deduce that $\Psi$ is a homomorphism.
\end{proof}

Now we are ready to show that the  homotopy $\huaO$-operators of weight $1$ induce \opt homotopy post-Lie algebras.

\begin{thm}\label{homotopy-RB-homotopy-post-lie}
Let $T=\sum_{i=0}^{+\infty}T_i\in \Hom^0(S(\h),\g)$ be a homotopy $\huaO$-operator of weight $1$ on a \sgla~ $(\g,[\cdot,\cdot]_\g)$ with respect to an action  $\rho:\g\lon\Der(\h)$. Then $(\h,[\cdot,\cdot]_\h,\{\oprn_k\}_{k=1}^\infty)$ is an \opt homotopy post-Lie algebra, where  $\oprn_k:\otimes^k \h\lon \h$ $(k\ge 1)$  are   linear maps  of degree $1$ defined by
\begin{eqnarray}
\label{homotopy-o-to-homotopy-post-lie}\oprn_k(v_1,\cdots,v_k):=\rho(T_{k-1}(v_1,\cdots,v_{k-1})\big)v_k,\quad \forall v_1\cdots,v_k\in\h.
\end{eqnarray}
\end{thm}

\begin{proof}
By Theorem \ref{hmotopy-o-operator-dgla} and Theorem \ref{homo-dg-lie}, we deduce that $\Psi(T)$ ia a Maurer-Cartan element of the \dgla~ $(s\bar{\huaC}^*(\h,\h),[\cdot,\cdot]^c,\partial)$. Moreover, by Theorem \ref{homotopy-post-dg-lie}, we obtain that $(\h,[\cdot,\cdot]_\h,\{\oprn_k\}_{k=1}^\infty)$ is an \opt homotopy post-Lie algebra.
\end{proof}

\begin{cor}\label{co:oprelie}
Let $T=\sum_{i=0}^{+\infty}T_i\in \Hom^0(S(V),\g)$ be a homotopy $\huaO$-operator of weight $0$ on a \sgla~ $(\g,[\cdot,\cdot]_\g)$ with respect to an action  $\rho:\g\lon\gl(V)$. Then $(V, \{\oprn_k\}_{k=1}^\infty)$ is a pre-Lie$_\infty$ algebra, where  $\oprn_k:\otimes^k V\lon V$ $(k\ge 1)$  are   linear maps  of degree $1$ defined by
\begin{eqnarray}
 \oprn_k(v_1,\cdots,v_k):=\rho(T_{k-1}(v_1,\cdots,v_{k-1})\big)v_k,\quad \forall v_1\cdots,v_k\in V.
\end{eqnarray}
\end{cor}
 It is straightforward to obtain the following result.

\begin{pro}\label{homotopy-post-lie-to-homotopy-lie}
  Let $(\g,[\cdot,\cdot]_\g,\{\oprn_k\}_{k=1}^\infty)$ be an \opt homotopy post-Lie algebra. Then $(\g,\{l_k\}_{k=1}^{\infty})$ is an $L_\infty$-algebra, where $l_1=\oprn_1$, $l_2$ is defined by
  \begin{equation}
    l_2(x,y)=\oprn_2(x,y)+(-1)^{xy}\oprn_2(y,x)+[x,y]_\g,
  \end{equation}
and for $k\geq3$, $l_k$ is defined by
\begin{equation}
l_k(x_1,\cdots,x_k)=\sum_{i=1}^{k}(-1)^{x_i(x_{i+1}+\cdots+x_k)}\oprn_k(x_{1},\cdots,\hat{x}_{i},\cdots,x_{k},x_i).
\end{equation}
\end{pro}

\begin{defi}{\rm (\cite{KS,MZ})}
Let $(V,\{l_k\}_{k=1}^\infty)$ be an $L_\infty$-algebra and $(V',\{l_m'\}_{m=1}^\infty)$ an $L_\infty$-algebra in which $l_m'=0$ for all $m\ge 1$ except $m=2$.
A curved $L_\infty$-algebra {\bf homomorphism} from $(V,\{l_k\}_{k=1}^\infty)$ to $(V',l_2')$ consists of a collection of degree $0$ graded multilinear maps $f_k:V^{\otimes k}\lon V',~ k\ge 0$ with the property that,
$f_n(v_{\sigma(1)},\cdots,v_{\sigma(n)})
=\varepsilon(\sigma)f_n(v_1,\cdots,v_n),
$ for any $n\geq 0$ and homogeneous elements $v_1,\cdots,v_n\in V$,
and
\begin{eqnarray*}
&&\sum_{i=1}^n\sum_{\sigma\in \mathbb S_{(i,n-i)} }\varepsilon(\sigma)f_{n-i+1}(l_i(v_{\sigma(1)},\cdots,v_{\sigma(i)}),v_{\sigma(i+1)},\cdots,v_{\sigma(n)})\\
&=&\frac{1}{2}\sum_{i=0}^n\sum_{\sigma\in \mathbb S_{(i,n-i)} }\varepsilon(\sigma)l_2'(f_{i}(v_{\sigma(1)},\cdots,v_{\sigma(i)}),f_{n-i}(v_{\sigma(i+1)},\cdots,v_{\sigma(n)})).
\end{eqnarray*}
\end{defi}

Then combining Theorem \ref{homotopy-RB-homotopy-post-lie} and Theorem \ref{homotopy-post-lie-to-homotopy-lie}, we obtain
\begin{thm}
Let $T=\sum_{i=0}^{+\infty}T_i\in \Hom(S(\h),\g)$ be a homotopy $\huaO$-operator of weight $1$ on a \sgla~ $(\g,[\cdot,\cdot]_\g)$ with respect to the action $\rho:\g\lon\Der(\h)$. Then $T$ is a curved $L_\infty$-algebra homomorphism from the  $L_\infty$-algebra $(\h,\{l_k\}_{k=1}^\infty)$ to $(\g,[\cdot,\cdot]_\g)$.
\end{thm}

\begin{proof}
For all $v_1,\cdots,v_n\in \h$, we have
\begin{eqnarray*}
&&\sum_{i=1}^n\sum_{\sigma\in \mathbb S_{(i,n-i)} }\varepsilon(\sigma)T_{n-i+1}(l_i(v_{\sigma(1)},\cdots,v_{\sigma(i)}),v_{\sigma(i+1)},\cdots,v_{\sigma(n)})\\
&=&\sum_{\sigma\in \mathbb S_{(1,n-1)} }\varepsilon(\sigma)T_{n}\Big(\rho(T_{0})v_{\sigma(1)},v_{\sigma(2)},\cdots,v_{\sigma(n)}\Big)\\
&&+\sum_{\sigma\in \mathbb S_{(2,n-1)}}\varepsilon(\sigma)T_{n-1}(\rho(T_1(v_{\sigma(1)}))v_{\sigma(2)}+(-1)^{v_{\sigma(1)}v_{\sigma(2)}}\rho(T_1(v_{\sigma(2)}))v_{\sigma(1)}+[v_{\sigma(1)},v_{\sigma(2)}]_\g,\cdots,v_{\sigma(n)})\\
&&+\sum_{i=3}^n\sum_{\sigma\in \mathbb S_{(i,n-i)} }\sum_{l=1}^{i}\varepsilon(\sigma)(-1)^{v_{\sigma(l)}(v_{\sigma(1+1)}+\cdots+v_{\sigma(i)})}T_{n-i+1}\Big(\rho(T_{i-1}(v_{\sigma(1)},\cdots,\hat{v}_{\sigma(l)},\cdots,v_{\sigma(i)})\big)\hat{v}_{\sigma(l)},v_{\sigma(i+1)},\cdots,v_{\sigma(n)}\Big)\\
&=&\sum_{1\le i< j\le n}(-1)^{v_i(v_1+\cdots+v_{i-1})+v_j(v_1+\cdots+v_{j-1})+v_iv_j}T_{n-1}([v_i,v_j]_{\h},v_1,\cdots,\hat{v}_i,\cdots,\hat{v}_{j},\cdots,v_{n})\\
&&\sum_{i=1}^n\sum_{\sigma\in \mathbb S_{(i-1,1,n-i)} }\varepsilon(\sigma)T_{n-i+1}\Big(\rho(T_{i-1}(v_{\sigma(1)},\cdots,v_{\sigma(i-1)})\big)v_{\sigma(i)},v_{\sigma(i+1)},\cdots,v_{\sigma(n)}\Big)\\
&\stackrel{\eqref{homotopy-rota-baxter-o}}{=}&\frac{1}{2}\sum_{i=0}^n\sum_{\sigma\in \mathbb S_{(i,n-i)} }\varepsilon(\sigma)[T_{i}(v_{\sigma(1)},\cdots,v_{\sigma(i)}),T_{n-i}(v_{\sigma(i+1)},\cdots,v_{\sigma(n)})]_\g,
\end{eqnarray*}
which implies that $T$ is a curved $L_\infty$-algebra homomorphism.
\end{proof}

Similarly, the above result also holds for homotopy $\huaO$-operators of weight 0.

\section{Classification of 2-term skeletal \opt homotopy post-Lie algebras}
\label{sec:repcoh}
In general, it is expected that the 2-term homotopy of an algebra structure is equivalent to the categorification of this algebraic structure, and  the 2-term homotopy algebras are quasi-isomorphic to the 2-term skeletal homotopy algebras, which are classified by the third cohomological group.
Baez and Crans~\cite{BC} accomplished these for Lie algebras. In this spirit, we classify  2-term skeletal \opt homotopy post-Lie algebras by the third cohomology group of a post-Lie algebra. For this purpose, we first define representations of post-Lie algebras and then develop the corresponding cohomology theory.

\subsection{Representations of post-Lie algebras}

Here we introduce the notion of a representation of a post-Lie algebra $(\g,[\cdot,\cdot]_\g,\rhd)$ on a vector space $V$. We show that there is naturally an induced representation of the subadjacent Lie algebra $\g^C$ on $\Der(\g,V)$. This fact plays a crucial role in our study of cohomology groups of post-Lie algebras in the next subsection.

\begin{defi}\label{defi:hom-pre representation}
 A {\bf representation} of a post-Lie algebra $(\g,[\cdot,\cdot]_\g,\rhd)$ on a vector space $V$ is a triple $(\rho,\mu,\nu)$, where $\rho:\g\longrightarrow \gl(V)$ is a representation of the Lie algebra $(\g,[\cdot,\cdot]_\g)$ on $V$, and $\mu,\nu:\g\longrightarrow \gl(V)$ are  linear maps satisfying
\begin{eqnarray}
\label{rep-1}\rho(x\rhd y)&=&\mu(x)\circ\rho(y)-\rho(y)\circ\mu(x),\\
\label{rep-2}\nu([x,y]_\g)&=&\rho(x)\circ\nu(y)-\rho(y)\circ\nu(x),\\
\label{rep-3}\mu([x,y]_\g)&=&\mu(x)\circ\mu(y)-\mu(x\rhd y)-\mu(y)\circ\mu(x)+\mu(y\rhd x),\\
\label{rep-4}\nu(y)\circ\rho(x)&=&\mu(x)\circ\nu(y)-\nu(y)\circ\mu(x)-\nu(x\rhd y)+\nu(y)\circ\nu(x), \quad \forall x,y\in \g.
\end{eqnarray}
\end{defi}

Let $(\g,[\cdot,\cdot]_\g,\rhd)$ be a post-Lie algebra and $(V;\rho,\mu,\nu)$ a representation of $(\g,[\cdot,\cdot]_\g,\rhd)$.
By Proposition \ref{double-Lie} and \eqref{rep-3}, we deduce that $(V;\mu)$ is a representation of the sub-adjacent Lie algebra $(\g,[\cdot,\cdot]_C)$. It is obvious that $(\g;\ad,L_\rhd,R_\rhd)$ is a representation of a post-Lie algebra on itself, which is called the {\bf regular representation}.

Let $(V;\rho,\mu,\nu)$ be a representation of a post-Lie algebra $(\g,[\cdot,\cdot]_\g,\rhd)$. Define a Lie bracket $[\cdot,\cdot]_\rho:\otimes^2(\g\oplus V)\lon \g\oplus V$  by
\begin{eqnarray}\label{semidirect-post-Lie-bracket-1}
[x_1+v_1,x_2+v_2]_{\rho}:=[x_1,x_2]_\g+\rho(x_1)v_2-\rho(x_2)v_1,
\end{eqnarray}
and a bilinear operation $\rhd_{\mu,\nu}:\otimes^2(\g\oplus V)\lon \g\oplus V$  by
\begin{eqnarray}\label{semidirect-post-Lie-bracket-2}
(x_1+v_1)\rhd_{\mu,\nu}(x_2+v_2):= x_1\rhd x_2+\mu(x_1)v_2+\nu(x_2)v_1.
\end{eqnarray}

By straightforward computations, we have
\begin{thm}\label{semidirect-post-Lie}
With the above notations, $(\g\oplus V,[\cdot,\cdot]_\rho,\rhd_{\mu,\nu})$ is a post-Lie algebra.
\end{thm}
This post-Lie algebra is called the
{\bf semidirect product} of the post-Lie algebra $(\g,[\cdot,\cdot]_\g,\rhd)$ and the representation $(V;\rho,\mu,\nu)$, and denoted by $\g\ltimes_{\rho,\mu,\nu}V$.

\begin{pro}\label{pro:representa}
Let $(V;\rho,\mu,\nu)$ be a representation of a post-Lie algebra $(\g,[\cdot,\cdot]_\g,\rhd)$. Then $(V;\rho+\mu-\nu)$ is a representation of the sub-adjacent Lie algebra $(\g,[\cdot,\cdot]_C)$.
\end{pro}

\begin{proof}
By Theorem \ref{semidirect-post-Lie}, we have the semidirect product post-Lie algebra $\g\ltimes_{\rho,\mu,\nu}V$. Considering its sub-adjacent Lie algebra structure $[\cdot,\cdot]_C$, we have
\begin{eqnarray}
\nonumber[(x_1+v_1),(x_2+v_2)]_C&=&(x_1+v_1)\rhd_{\mu,\nu}(x_2+v_2)-(x_2+v_2)\rhd_{\mu,\nu}(x_1+v_1)+[x_1+v_1,x_2+v_2]_{\rho}\\
                               \nonumber&=&x_1\rhd x_2+\mu(x_1)v_2+\nu(x_2)v_1-x_2\rhd x_1-\mu(x_2)v_1-\nu(x_1)v_2\\
                               \nonumber&&+[x_1,x_2]_\g+\rho(x_1)v_2-\rho(x_2)v_1\\
                               &=&[x_1,x_2]_C+(\rho+\mu-\nu)(x_1)v_2-(\rho+\mu-\nu)(x_2)v_1.\label{eq:samesubadj}
\end{eqnarray}
Thus,  $(V;\rho+\mu-\nu)$ is a representation of the sub-adjacent Lie algebra $(\g,[\cdot,\cdot]_C)$.
\end{proof}

If $(\rho,\mu,\nu)=(\ad,L_\rhd,R_\rhd)$ is the regular representation of a post-Lie algebra $(\g,[\cdot,\cdot]_\g,\rhd)$, then $\ad+L_\rhd-R_\rhd$ is the adjoint representation of the sub-adjacent Lie algebra $(\g,[\cdot,\cdot]_C)$.

\begin{cor}\label{sub-adjacent-3-Lie}
Let $(V;\rho,\mu,\nu)$ be a representation of a post-Lie algebra $(\g,[\cdot,\cdot]_\g,\rhd)$. Then the semidirect product post-Lie algebras $\g\ltimes_{\rho,\mu,\nu}V$ and $\g\ltimes_{0,\rho+\mu-\nu,0}V$ given by the representations $(V;\rho,\mu,\nu)$ and $(V;0,\rho+\mu-\nu,0)$ respectively have the same sub-adjacent Lie algebra $\g^C\ltimes_{\rho+\mu-\nu}V$ given by \eqref{eq:samesubadj}, which
  is the semidirect product of the Lie algebra $(\g,[\cdot,\cdot]_C)$ and its representation   $(V;\rho+\mu-\nu)$.
\end{cor}

Let $(V;\rho,\mu,\nu)$ be a representation of a post-Lie algebra $(\g,[\cdot,\cdot]_\g,\rhd)$. We set
\begin{eqnarray}\label{derivation-rep}
\Der(\g,V):=\Big\{f\in\Hom(\g,V)\big|f([x,y]_\g)=\rho(x)f(y)-\rho(y)f(x),~\forall x,y\in\g\Big\}
\end{eqnarray}
and define $\hat{\rho}:\g\lon\Hom\Big(\Der(\g,V),\Hom(\g,V)\Big)$ by
\begin{eqnarray}\label{important-rep}
\Big(\hat{\rho}(x)(f)\Big)y:=\mu(x)f(y)+\nu(y)f(x)-f(x\rhd y),\quad\forall x,y\in\g,~f\in\Der(\g,V).
\end{eqnarray}

By a straightforward computation, we deduce
\begin{lem}
 For all $x\in\g$, we have $\hat{\rho}(x)\in \gl(\Der(\g,V))$.
\end{lem}

\begin{thm}\label{thm:rep}
Let $(V;\rho,\mu,\nu)$ be a representation of a post-Lie algebra $(\g,[\cdot,\cdot]_\g,\rhd)$. Then $(\Der(\g,V);\hat{\rho})$ is a representation of the sub-adjacent Lie algebra $(\g,[\cdot,\cdot]_C)$, where $\hat{\rho}$ is given by \eqref{important-rep}.
\end{thm}

\begin{proof}
By \eqref{important-rep}, for all $x,y,z\in\g$ and $f\in\Der(\g,V)$, we have
\begin{eqnarray*}
&&\Big(\big([\hat{\rho}(x),\hat{\rho}(y)]-\hat{\rho}([x,y]_C)\big)(f)\Big)z\\
&=&\hat{\rho}(x)(\hat{\rho}(y)(f))z-\hat{\rho}(y)(\hat{\rho}(x)(f))z-\hat{\rho}([x,y]_C)(f)z\\
&=&\mu(x)(\hat{\rho}(y)(f))z+\nu(z)(\hat{\rho}(y)(f))x-(\hat{\rho}(y)(f))(x\rhd z)\\
&&-\mu(y)(\hat{\rho}(x)(f))z-\nu(z)(\hat{\rho}(x)(f))y+(\hat{\rho}(x)(f))(y\rhd z)\\
&&-\mu([x,y]_C)f(z)-\nu(z)f([x,y]_C)+f([x,y]_C\rhd z)\\
&=&\mu(x)\Big(\mu(y)f(z)+\nu(z)f(y)-f(y\rhd z)\Big)+\nu(z)\Big(\mu(y)f(x)+\nu(x)f(y)-f(y\rhd x)\Big)\\
&&-\Big(\mu(y)f(x\rhd z)+\nu(x\rhd z)f(y)-f(y\rhd(x\rhd z))\Big)-\mu(y)\Big(\mu(x)f(z)+\nu(z)f(x)-f(x\rhd z)\Big)\\
&&-\nu(z)\Big(\mu(x)f(y)+\nu(y)f(x)-f(x\rhd y)\Big)+\Big(\mu(x)f(y\rhd z)+\nu(y\rhd z)f(x)-f(x\rhd(y\rhd z))\Big)\\
&&-\mu\Big(x\rhd y-y\rhd x+[x,y]_\g\Big)f(z)-\nu(z)f(x\rhd y-y\rhd x+[x,y]_\g)\\
&&+f\Big((x\rhd y-y\rhd x+[x,y]_\g)\rhd z\Big)\\
&\stackrel{\eqref{Post-2},\eqref{rep-3}}{=}&\mu(x)\nu(z)f(y)+\nu(z)\mu(y)f(x)+\nu(z)\nu(x)f(y)-\nu(x\rhd z)f(y)\\
&&-\mu(y)\nu(z)f(x)-\nu(z)\mu(x)f(y)-\nu(z)\nu(y)f(x)+\nu(y\rhd z)f(x)-\nu(z)f([x,y]_\g)\\
&\stackrel{\eqref{derivation-rep}}{=}&\underline{\mu(x)\nu(z)f(y)}\underbrace{+\nu(z)\mu(y)f(x)}\underline{+\nu(z)\nu(x)f(y)-\nu(x\rhd z)f(y)}\\
&&\underbrace{-\mu(y)\nu(z)f(x)}\underline{-\nu(z)\mu(x)f(y)}\underbrace{-\nu(z)\nu(y)f(x)+\nu(y\rhd z)f(x)}\underline{-\nu(z)\rho(x)f(y)}\underbrace{+\nu(z)\rho(y)f(x)}\\
&\stackrel{\eqref{rep-4}}{=}&0.
\end{eqnarray*}
Thus,   $(\Der(\g,V);\hat{\rho})$ is a representation of the sub-adjacent Lie algebra $(\g,[\cdot,\cdot]_C)$.
\end{proof}

\subsection{Cohomology groups of post-Lie algebras} \label{subsec:cohomologyP}
In this subsection, we define the cohomology groups of a post-Lie algebra with coefficients in an arbitrary representation. Furthermore, we establish precise relationship between the cohomology groups of a  post-Lie algebra and those of its subadjacent Lie algebra.

Let $(V;\rho,\mu,\nu)$ be a representation of a post-Lie algebra $(\g,[\cdot,\cdot]_\g,\rhd)$. We have the natural isomorphism $$\Phi:\Hom(\wedge^{n-1} \g\otimes \g,V)\longrightarrow \Hom(\wedge^{n-1} \g, \Hom(\g,V))$$ defined by
\begin{equation}\label{eq:3.1}
\Big(\Phi(\omega)(x_1,\cdots,x_{n-1})\Big)x_n:=\omega(x_1,\cdots,x_{n-1},x_n),\quad \forall x_1,\cdots,x_{n-1},x_n\in \g.
\end{equation}
Define the set of $0$-cochains to be $0$. For $n\ge1$, we define the set of $n$-cochains $C_{\Der}^{n}(\g,V)$ by
\begin{equation}
 C_{\Der}^{n}(\g,V):=\Big\{f\in\Hom(\wedge^{n-1} \g\otimes \g,V)\Big|\Phi_{}(f)(x_1,\cdots,x_{n-1})\in\Der(\g,V),\forall x_1,\cdots,x_{n-1}\in \g\Big\}.
\end{equation}
For all $f \in C_{\Der}^n(\g,V),~ x_1,\cdots,x_{n+1} \in \g$, define the operator
$\delta:C_{\Der}^n(\g,V)\longrightarrow \Hom(\wedge^{n} \g\otimes \g,V)$ by
\begin{eqnarray}\notag
&&(\delta f)(x_1,\cdots,x_{n+1})=\sum_{i=1}^{n}(-1)^{i+1}\mu(x_i)f(x_1,\cdots,\hat{x}_i,\cdots, x_{n},x_{n+1})\\
\label{eq:12.1}&&\qquad+\sum_{i=1}^{n}(-1)^{i+1}\nu(x_{n+1})f(x_1,\cdots,\hat{x}_i,\cdots, x_{n},x_i)-\sum_{i=1}^{n}(-1)^{i+1}f(x_1,\cdots,\hat{x}_i,\cdots, x_{n},x_i\rhd x_{n+1})\\
  \nonumber&&\qquad+\sum_{1\le i<j\le n}(-1)^{i+j}f(x_i\rhd x_j-x_j\rhd x_i+[x_i,x_j]_\g,x_1,\cdots,\hat{x}_i,\cdots,\hat{x}_j,\cdots, x_{n},x_{n+1}).
\end{eqnarray}

\begin{pro}
  For all $f\in C_{\Der}^{n}(\g,V)$, we have $\delta f\in C_{\Der}^{n+1}(\g,V)$.
\end{pro}

\begin{proof}
By the definition of $C_{\Der}^{n+1}(\g,V)$, we just need to prove that
$
\Phi(\delta f)(x_1,\cdots,x_{n})$ is in $\Der(\g,V)$ for any $x_1,\cdots,x_{n}\in \g.$
For all $x,y\in\g$, we have
\begin{eqnarray*}
\lefteqn{\Phi(\delta f)(x_1,\cdots,x_{n})([x,y]_\g)
\stackrel{\eqref{eq:3.1}}{=}(\delta f)(x_1,\cdots,x_{n},[x,y]_\g)}\\
&\stackrel{\eqref{eq:12.1}}{=}&\sum_{i=1}^{n}(-1)^{i+1}\mu(x_i)f(x_1,\cdots,\hat{x}_i,\cdots, x_{n},[x,y]_\g) +\sum_{i=1}^{n}(-1)^{i+1}\nu([x,y]_\g)f(x_1,\cdots,\hat{x}_i,\cdots, x_{n},x_i)\\
&&-\sum_{i=1}^{n}(-1)^{i+1}f(x_1,\cdots,\hat{x}_i,\cdots, x_{n},x_i\rhd [x,y]_\g)\\
&&+\sum_{1\le i<j\le n}(-1)^{i+j}f(x_i\rhd x_j-x_j\rhd x_i+[x_i,x_j]_\g,x_1,\cdots,\hat{x}_i,\cdots,\hat{x}_j,\cdots, x_{n},[x,y]_\g)\\
&\stackrel{\eqref{Post-1},\eqref{rep-2}}{=}&\sum_{i=1}^{n}(-1)^{i+1}\mu(x_i)\Big(\rho(x)f(x_1,\cdots,\hat{x}_i,\cdots, x_{n},y)-\rho(y)f(x_1,\cdots,\hat{x}_i,\cdots, x_{n},x)\Big)\\
&&+\sum_{i=1}^{n}(-1)^{i+1}\Big(\rho(x)\nu(y)-\rho(y)\nu(x)\Big)f(x_1,\cdots,\hat{x}_i,\cdots, x_{n},x_i)\\
&& -\sum_{i=1}^{n}(-1)^{i+1}\rho(x_i\rhd x)f(x_1,\cdots,\hat{x}_i,\cdots, x_{n},y)+\sum_{i=1}^{n}(-1)^{i+1}\rho(y)f(x_1,\cdots,\hat{x}_i,\cdots, x_{n},x_i\rhd x)\\
&&-\sum_{i=1}^{n}(-1)^{i+1}\rho(x)f(x_1,\cdots,\hat{x}_i,\cdots, x_{n},x_i\rhd y) +\sum_{i=1}^{n}(-1)^{i+1}\rho(x_i\rhd y)f(x_1,\cdots,\hat{x}_i,\cdots, x_{n},x) \\
&&+\sum_{1\le i<j\le n}(-1)^{i+j}\rho(x)f(x_i\rhd x_j-x_j\rhd x_i+[x_i,x_j]_\g,x_1,\cdots,\hat{x}_i,\cdots,\hat{x}_j,\cdots, x_{n},y)\\
&&-\sum_{1\le i<j\le n}(-1)^{i+j}\rho(y)f(x_i\rhd x_j-x_j\rhd x_i+[x_i,x_j]_\g,x_1,\cdots,\hat{x}_i,\cdots,\hat{x}_j,\cdots, x_{n},x)\\
&\stackrel{\eqref{rep-1}}{=}&\sum_{i=1}^{n}(-1)^{i+1}\Big(\rho(x)\nu(y)-\rho(y)\nu(x)\Big)f(x_1,\cdots,\hat{x}_i,\cdots, x_{n},x_i)\\
&&+\sum_{i=1}^{n}(-1)^{i+1}\rho(x)\Big(\mu(x_i)f(x_1,\cdots,\hat{x}_i,\cdots, x_{n},y)\Big)+\sum_{i=1}^{n}(-1)^{i+1}\rho(y)f(x_1,\cdots,\hat{x}_i,\cdots, x_{n},x_i\rhd x)\\
&&-\sum_{i=1}^{n}(-1)^{i+1}\rho(x)f(x_1,\cdots,\hat{x}_i,\cdots, x_{n},x_i\rhd y)-\sum_{i=1}^{n}(-1)^{i+1}\rho( y)\Big(\mu(x_i)f(x_1,\cdots,\hat{x}_i,\cdots, x_{n},x)\Big)\\
&&+\sum_{1\le i<j\le n}(-1)^{i+j}\rho(x)f(x_i\rhd x_j-x_j\rhd x_i+[x_i,x_j]_\g,x_1,\cdots,\hat{x}_i,\cdots,\hat{x}_j,\cdots, x_{n},y)\\
&&-\sum_{1\le i<j\le n}(-1)^{i+j}\rho(y)f(x_i\rhd x_j-x_j\rhd x_i+[x_i,x_j]_\g,x_1,\cdots,\hat{x}_i,\cdots,\hat{x}_j,\cdots, x_{n},x)\\
&=&\rho(x)\Big(\Phi(\delta f)(x_1,\cdots,x_{n})y\Big)-\rho(y)\Big(\Phi(\delta f)(x_1,\cdots,x_{n})x\Big).
\end{eqnarray*}
Thus we deduce that $\Phi(\delta f)(x_1,\cdots,x_{n})$ is in $\Der(\g,V)$.
\end{proof}

To prove that the operator $\delta$ is indeed a coboundary operator, i.e. $\delta\circ \delta=0,$ we need some preparations.
\begin{pro}\label{pro:commutative}
Let $(V;\rho,\mu,\nu)$ be a representation of a post-Lie algebra $(\g,[\cdot,\cdot]_\g,\rhd)$.
Then we have $\Phi \circ \delta=d_{\hat{\rho}} \circ \Phi$,
where $d_{\hat{\rho}}$ is the coboundary operator of the sub-adjacent Lie algebra $\g^C$ with coefficients in the representation $(\Der(\g,V);\hat{\rho})$ given in Theorem \ref{thm:rep} and $\delta$ is defined by \eqref{eq:12.1}.
\end{pro}

\begin{proof}
For all  $f\in C_{\Der}^{n}(\g,V)$ and $x_1,\cdots,x_{n+1}\in\g$, we have
\begin{eqnarray*}
&&\Big(d_{\hat{\rho}}(\Phi_{}(f))(x_1,\cdots,x_{n})\Big)x_{n+1}\\
&=&\sum_{i=1}^{n}(-1)^{i+1}\Big(\hat{\rho}(x_i)\Phi_{}(f)(x_1,\cdots,\hat{x}_i,\cdots,x_n)\Big)x_{n+1}\\
&&+\sum_{1\le i<j\le n}(-1)^{i+j}\Big(\Phi_{}(f)([x_i,x_j]_C,x_1,\cdots,\hat{x}_i,\cdots,\hat{x}_j,\cdots,x_n)\Big)x_{n+1}\\
&\stackrel{\eqref{important-rep}}{=}&\sum_{i=1}^{n}(-1)^{i+1}\mu(x_i)\Phi_{}(f)(x_1,\cdots,\hat{x}_i,\cdots,x_n)x_{n+1}+\sum_{i=1}^{n}(-1)^{i+1}\nu(x_{n+1})\Phi_{}(f)(x_1,\cdots,\hat{x}_i,\cdots,x_n)x_{i}\\
&&-\sum_{i=1}^{n}(-1)^{i+1}\Phi_{}(f)(x_1,\cdots,\hat{x}_i,\cdots,x_n)(x_i\rhd x_{n+1})\\
&&+\sum_{1\le i<j\le n}(-1)^{i+j}\Big(\Phi_{}(f)(x_i\rhd x_j-x_j\rhd x_i+[x_i,x_j]_\g,x_1,\cdots,\hat{x}_i,\cdots,\hat{x}_j,\cdots,x_n)\Big)x_{n+1}\\
&\stackrel{\eqref{eq:3.1}}{=}&\sum_{i=1}^{n}(-1)^{i+1}\mu(x_i)f(x_1,\cdots,\hat{x}_i,\cdots,x_n,x_{n+1})+\sum_{i=1}^{n}(-1)^{i+1}\nu(x_{n+1})f(x_1,\cdots,\hat{x}_i,\cdots,x_n,x_{i})\\
&&-\sum_{i=1}^{n}(-1)^{i+1}f(x_1,\cdots,\hat{x}_i,\cdots,x_n,x_i\rhd x_{n+1})\\
&&+\sum_{1\le i<j\le n}(-1)^{i+j}f(x_i\rhd x_j-x_j\rhd x_i+[x_i,x_j]_\g,x_1,\cdots,\hat{x}_i,\cdots,\hat{x}_j,\cdots,x_n,x_{n+1})\\
&\stackrel{\eqref{eq:12.1}}{=}&(\delta f)(x_1,\cdots,x_{n+1})\\
&\stackrel{\eqref{eq:3.1}}{=}&\Big(\Phi_{}(\delta f)(x_1,\cdots,x_{n})\Big)x_{n+1},
\end{eqnarray*}
which implies that $d_{\hat{\rho}} \circ \Phi=\Phi_{}\circ \delta $. \end{proof}

\begin{thm}\label{cohomology-post-lie}
The operator $\delta:C_{\Der}^n(\g,V)\longrightarrow C_{\Der}^{n+1}(\g,V)$ defined by \eqref{eq:12.1} satisfies $\delta\circ\delta=0$.
\end{thm}
\begin{proof}
By Proposition \ref{pro:commutative}, we have $\delta=\Phi^{-1}\circ d_{\hat{\rho}} \circ \Phi$.
Thus, by the fact that $d_{\hat{\rho}} \circ d_{\hat{\rho}}=0$, we obtain $\delta\circ\delta=\Phi^{-1}\circ d_{\hat{\rho}} \circ d_{\hat{\rho}}\circ\Phi=0.$
\end{proof}

Let $(V;\rho,\mu,\nu)$ be a representation of a post-Lie algebra $(\g,[\cdot,\cdot]_\g,\rhd)$. Denote by $C_{\Der}^*(\g,V):=\bigoplus_{n\geq 1} C_{\Der}^n(\g,V)$. Then we have the cochain complex $(C_{\Der}^*(\g,V),\delta)$. Denote the set of closed $n$-cochains by $Z^n(\g,V)$ and the set of exact $n$-cochains by $B^n(\g,V)$. We denote by $H^n(\g,V)=Z^n(\g,V)/B^n(\g,V)$, and call them the  {\bf cohomology groups of the post-Lie algebra} $(\g,[\cdot,\cdot]_\g,\rhd)$ with coefficients in the representation $(V;\rho,\mu,\nu)$.

It is obvious that $f\in C_{\Der}^{1}(\g,V)$ is closed if and only if $f\in\Der(\g,V)$ and
         $$
          \nu(y)f(x)+\mu(x)f(y)-f(x\rhd y)=0,\quad\forall x,y\in\g.
         $$
Also  $f\in C_{\Der}^2(\g,V)$ is closed if and only if $\Phi(f)\in \Hom(\g,\Der(\g,V))$ and
\begin{eqnarray*}
&&\nu(x_3)f(x_2,x_1)+\mu(x_1)f(x_2,x_3)-f(x_2,x_1\rhd x_3)-\nu(x_3)f(x_1,x_2)-\mu(x_2)f(x_1,x_3)\\&&+f(x_1,x_2\rhd x_3)
-f(x_1\rhd x_2,x_3)+f(x_2\rhd x_1,x_3)-f([x_1,x_2]_\g,x_3)=0,\quad\forall x_1,x_2,x_3\in\g.
\end{eqnarray*}

There is a close relationship between the cohomology groups of post-Lie algebras and those of the corresponding sub-adjacent Lie algebras.

\begin{thm}
Let $(V;\rho,\mu,\nu)$ be a representation of a post-Lie algebra $(\g,[\cdot,\cdot]_\g,\rhd)$. Then the cohomology group $H^n(\g,V)$ of the post-Lie algebra $(\g,[\cdot,\cdot]_\g,\rhd)$  and the cohomology group $H^{n-1}(\g^C,\Der(\g,V))$ of the subadjacent Lie algebra $\g^C$ are isomorphic for all $n\geq 1$.
\end{thm}
\begin{proof}
By Proposition \ref{pro:commutative}, we deduce that $\Phi$ is an isomorphism from the cochain complex $\big(C_{\Der}^\ast(\g,V),\delta\big)$ to the cochain complex $\big(\huaC^{\ast-1}(\g,\Der(\g,V)),d_{\hat{\rho}}\big)$. Thus, $\Phi$ induces an isomorphism $\Phi_\ast$ from $H^\ast(\g,V)$ to $H^{\ast-1}(\g^C,\Der(\g,V))$.
\end{proof}

In the above theorem, if $[\cdot,\cdot]_\frakg$ and $\rho$ are zero, then the post-Lie algebra is a pre-Lie algebra and we recover the result of~\cite{DA} as follows.

\begin{cor}
Let $(V;\mu,\nu)$ be a representation of a pre-Lie algebra $(\g,\rhd)$. Then the cohomology group $H^n(\g,V)$ of the pre-Lie algebra $(\g,\rhd)$  and the cohomology group $H^{n-1}(\g^C,\Hom(\g,V))$ of the subadjacent Lie algebra $\g^C$ are isomorphic for all $n\geq 1$.
\end{cor}

\subsection{Classification of 2-term skeletal \opt homotopy post-Lie algebras}
\label{ss:deform}

In this subsection, we first give an equivalent definition of an \opt homotopy post-Lie algebra and then classify  2-term skeletal \opt homotopy post-Lie algebras using the third cohomology group given in Section~\ref{subsec:cohomologyP}

For all $i\ge 1$, let $\Oprn_i:\wedge^{i-1}\g\otimes\g\lon\g$ be a graded linear map of degree $2-i$. Define $D(\Oprn_i):\odot^{i-1}s^{-1}\g\otimes s^{-1}\g\lon s^{-1}\g$ by
\begin{eqnarray*}
D(\Oprn_i)=(-1)^{\frac{i(i-1)}{2}}s^{-1}\circ \Oprn_i\circ s^{\otimes i},
\end{eqnarray*}
which is a graded linear map of degree $1$. This is illustrated by the following commutative diagram:
\[\begin{CD}
\wedge^{i-1}\g\otimes\g @>\Oprn_i>>                                    {\g}            \\
 @V\otimes^i s^{-1} VV                             @VV s^{-1} V \\
\odot^{i-1}s^{-1}\g\otimes s^{-1}\g@>D(\Oprn_i)>>                    {s^{-1}\g}
\end{CD}\]

Using this process, we can give an equivalent definition of an \opt homotopy post-Lie algebra as follows.
\begin{defi}
An {\bf \opt homotopy post-Lie algebra} is a graded Lie algebra $(\g,[\cdot,\cdot]_\g)$ equipped with a collection of linear maps $\Oprn_k:\otimes^k \g\lon \g, k\ge 1,$ of degree $2-k$ satisfying, for any homogeneous elements $x_1,\cdots,x_n,x_{n+1}\in \g$, the following conditions hold:
\begin{itemize}
\item[\rm(i)]
{\bf (graded antisymmetry)} for every $\sigma\in\mathbb S_{n-1},~n\ge1$,
\begin{eqnarray}
 \label{homotopy-post-lie-1}\Oprn_n(x_{\sigma(1)},\cdots,x_{\sigma(n-1)},x_n)=\chi(\sigma)\Oprn_n(x_1,\cdots,x_{n-1},x_n),
\end{eqnarray}
\item[\rm(ii)]
{\bf (graded derivation)} for all $n\ge 1$,
\begin{eqnarray}
 \Oprn_n(x_1,\cdots,x_{n-1},[x_n,x_{n+1}]_\g)
\label{homotopy-post-lie-2}&=&[\Oprn_n(x_1,\cdots,x_{n-1},x_n),x_{n+1}]_\g\\
 \nonumber&&+(-1)^{x_n(x_1+\cdots+x_{n-1}+n)}[x_n,\Oprn_n(x_1,\cdots,x_{n-1},x_{n+1})]_\g,
\end{eqnarray}
\item[\rm(iii)] for all $n\ge 1$,
\begin{eqnarray}
\nonumber&&\sum_{1\le i< j\le n-1}(-1)^{\beta}\Oprn_{n-1}([x_i,x_j]_\g,x_1,\cdots,\hat{x}_i,\cdots,\hat{x}_j,\cdots,x_n)\\
  \label{homotopy-post-lie1}&=&\sum_{i+j=n+1\atop i\ge1,j\geq2}\sum_{\sigma\in\mathbb S_{(i-1,1,j-2)}}(-1)^{i(j-1)}\chi(\sigma)\Oprn_j(\Oprn_i(x_{\sigma(1)},\cdots,x_{\sigma(i-1)},x_{\sigma(i)}), x_{\sigma(i+1)},\cdots,x_{\sigma(n-1)},x_{n})\\
\nonumber&&+\sum_{i+j=n+1\atop i\geq1,j\geq1}\sum_{\sigma\in\mathbb S_{(j-1,i-1)}}(-1)^{j-1}(-1)^{\alpha}\chi(\sigma)\Oprn_j(x_{\sigma(1)},\cdots,x_{\sigma(j-1)},\Oprn_i( x_{\sigma(j)},\cdots,x_{\sigma(n-1)},x_{n})),
\end{eqnarray}
where $\beta=x_i(x_1+\cdots+x_{i-1})+x_j(x_1+\cdots+x_{j-1})+x_ix_j+i+j$ and $\alpha=i(x_{\sigma(1)}+x_{\sigma(2)}+\cdots+x_{\sigma(j-1)})$.
\end{itemize}
\end{defi}

By \eqref{homotopy-post-lie-2}, for $n=1$, we have
\begin{eqnarray}\label{eq:coho1}
\Oprn_1([x_1,x_2]_\g)=[\Oprn_1(x_1),x_2]_\g+(-1)^{x_1}[x_1,\Oprn_1(x_2)]_\g.
\end{eqnarray}
By \eqref{homotopy-post-lie1}, for $n=1$, we have $\Oprn_1^2=0,$ which implies  that $(\g,\Oprn_1)$ is a complex. By \eqref{homotopy-post-lie1}, for $n=2$, we have
\begin{eqnarray}\label{eq:coho2}
0=-\Oprn_2(\Oprn_1(x_1),x_2)-(-1)^{x_1}\Oprn_2(x_1,\Oprn_1(x_2))+\Oprn_1(\Oprn_2(x_1,x_2)).
\end{eqnarray}

Now we will show that the corresponding cohomology space  $H^*(\g)$ of the complex $(\g,\Oprn_1)$  enjoys a graded post-Lie algebra structure and this justifies our definition of an ``\opt homotopy post-Lie algebra''.

\begin{thm}
Let $(\g,[\cdot,\cdot]_\g,\{\Oprn_k\}_{k=1}^{\infty})$ be an \opt homotopy post-Lie algebra. Then the cohomology space $H^*(\g)$ is a graded post-Lie algebra.
\end{thm}
\begin{proof}
 For any homogeneous element $x\in\ker(\Oprn_1)$, we denote by $\overline{x}\in H^*(\g)$ its cohomological class. First we define a graded bracket operation $[\cdot,\cdot]$  on the graded vector space $H^*(\g)$ by
\begin{eqnarray*}
[\bar{x},\bar{y}]:=\overline{[x,y]_\g},\quad\forall \bar{x},\bar{y}\in H^*(\g).
\end{eqnarray*}
If $\bar{x}=\bar{x'}$, then there exists $X\in\g$ such that $x'=x+\Oprn_1(X)$. Hence by \eqref{eq:coho1}, we have
$$
[\bar{x'},\bar{y}]=\overline{[x',y]_\g}=\overline{[x+\Oprn_1(X),y]_\g} =\overline{[x,y]_\g}+\overline{\Oprn_1([X,y]_\g)}=\overline{[x,y]_\g}=[\bar{x},\bar{y}],
$$
which implies that $[\cdot,\cdot]$ is well-defined. It is straightforward to obtain that $(H^*(\g),[\cdot,\cdot])$ is a graded Lie algebra.

Then we define a multiplication $\rhd$ on the graded vector space $H^*(\g)$ by
\begin{eqnarray*}
\bar{x}\rhd \bar{y}:=\overline{\Oprn_2(x,y)},\quad\forall \bar{x},\bar{y}\in H^*(\g).
\end{eqnarray*}
Similarly, by \eqref{eq:coho2}, we can deduce that $\rhd$ is well-defined.
By \eqref{homotopy-post-lie-2} for $n=2$, we have
\begin{eqnarray*}
 \bar{x}\rhd [\bar{y},\bar{z}]=\overline{\Oprn_2(x,[y,z]_\g)}=\overline{[\Oprn_2(x,y),z]_\g}+(-1)^{xy}\overline{[y,\Oprn_2(x,z)]_\g}= [\bar{x}\rhd \bar{y},\bar{z}]+ +(-1)^{xy}  [ \bar{y},\bar{x}\rhd\bar{z}].
\end{eqnarray*}
 Similarly,
by \eqref{homotopy-post-lie1} for $n=3$, we have
\begin{eqnarray*}
[\bar{x},\bar{y}]\rhd \bar{z}= a_\rhd(\bar{x},\bar{y},\bar{z})-a_\rhd(\bar{y},\bar{x},\bar{z}).
\end{eqnarray*}
Therefore,  $(H^*(\g),[\cdot,\cdot],\rhd)$ is a graded post-Lie algebra.
\end{proof}

By truncation, we obtain the definition of a 2-term \opt homotopy post-Lie algebra.

\begin{defi}\label{defi:2-term}
A {\bf 2-term \opt homotopy post-Lie algebra} is a $2$-term graded Lie algebra $(\g=\g_{0}\oplus \g_{-1},[\cdot,\cdot]_\g)$ equipped with
\begin{itemize}
  \item[$\bullet$] a linear map $\Oprn_1:\g_{-1}\lon \g_{0}$;
  \item[$\bullet$] a linear map $\Oprn_2:\g_{i}\otimes \g_{j}\lon \g_{i+j}$, $-1\leq i+j\leq 0$;
  \item[$\bullet$] a linear map $\Oprn_3:\wedge^2\g_{0}\otimes \g_{0}\lon \g_{-1}$
\end{itemize}
such that for all $x,y,z,w\in \g_{0}$ and $a,b\in\g_{-1}$, the following equalities hold:
\begin{itemize}
\item[\rm($a_1$)]  $\Oprn_1([x,a]_\g)=[x,\Oprn_1(a)]_\g$;
\item[\rm($a_2$)]$[\Oprn_1(a),b]_\g=[a,\Oprn_1(b)]_\g$;

\item[\rm($b_1$)]  $\Oprn_2(x,[y,z]_\g)=[\Oprn_2(x,y),z]_\g+[y,\Oprn_2(x,z)]_\g$;
\item[\rm($b_2$)]  $\Oprn_2(x,[y,a]_\g)=[\Oprn_2(x,y),a]_\g+[y,\Oprn_2(x,a)]_\g$;
\item[\rm($b_3$)]  $\Oprn_2(a,[x,y]_\g)=[\Oprn_2(a,x),y]_\g+[x,\Oprn_2(a,y)]_\g$;

\item[\rm($c$)]  $\Oprn_3(x,y,[z,w]_\g)=[\Oprn_3(x,y,z),w]_\g+[z,\Oprn_3(x,y,w)]_\g$;

\item[\rm($d_1$)] $\Oprn_1\Oprn_2(x,a)=\Oprn_2(x,\Oprn_1(a))$;

\item[\rm($d_2$)] $\Oprn_1\Oprn_2(a,x)=\Oprn_2(\Oprn_1(a),x)$;

\item[\rm($d_3$)] $\Oprn_2(\Oprn_1(a),b)=\Oprn_2(a,\Oprn_1(b))$;

\item[\rm($e_1$)] $\Oprn_2(x,\Oprn_2(y,z))-\Oprn_2(\Oprn_2(x,y),z)-\Oprn_2(y,\Oprn_2(x,z))+\Oprn_2(\Oprn_2(y,x),z)
    -\Oprn_2([x,y]_\g,z)=\Oprn_1\Oprn_3(x,y,z)$;

    \item[\rm($e_2$)] $\Oprn_2(x,\Oprn_2(y,a))-\Oprn_2(\Oprn_2(x,y),a)-\Oprn_2(y,\Oprn_2(x,a))+\Oprn_2(\Oprn_2(y,x),a)
            -\Oprn_2([x,y]_\g,a)= \Oprn_3(x,y,\Oprn_1(a))$;

  \item[\rm($e_3$)] $\Oprn_2(a,\Oprn_2(y,z))-\Oprn_2(\Oprn_2(a,y),z)-\Oprn_2(y,\Oprn_2(a,z))+\Oprn_2(\Oprn_2(y,a),z)
    -\Oprn_2([a,y]_\g,z)=\Oprn_3(\Oprn_1 (a),y,z)$;

      \item[\rm($f$)]     $\Oprn_2(x,\Oprn_3(y,z,w))-\Oprn_2(y,\Oprn_3(x,z,w))+\Oprn_2(z,\Oprn_3(x,y,w))
+\Oprn_2(\Oprn_3(y,z,x),w)-\Oprn_2(\Oprn_3(x,z,y),w)+\Oprn_2(\Oprn_3(x,y,z),w)
              -\Oprn_3(\Oprn_2 (x,y)-\Oprn_2 (y,x)+[x,y]_\g,z,w)
      -\Oprn_3(\Oprn_2 (y,z)-\Oprn_2 (z,y)+[y,z]_\g,x,w)
      +\Oprn_3(\Oprn_2 (x,z)-\Oprn_2 (z,x)+[x,z]_\g,y,w)
-\Oprn_3(y,z,\Oprn_2 (x,w))+\Oprn_3(x,z,\Oprn_2 (y,w))-\Oprn_3(x,y,\Oprn_2 (z,w))=0.
      $

    \end{itemize}
   A  $2$-term \opt homotopy post-Lie algebra $(\g=\g_{0}\oplus \g_{-1},[\cdot,\cdot]_\g,\Oprn_1,\Oprn_2,\Oprn_3)$ is said to be {\bf skeletal} if $\Oprn_1=0.$
\end{defi}

\begin{rmk}
If the underlying graded Lie algebra $(\g_{0}\oplus \g_{-1},[\cdot,\cdot]_\g)$ in a $2$-term \opt homotopy post-Lie algebra $(\g=\g_{0}\oplus \g_{-1},[\cdot,\cdot]_\g,\Oprn_1,\Oprn_2,\Oprn_3)$ is   abelian, then $(\g_{0}\oplus \g_{-1}, \Oprn_1,\Oprn_2,\Oprn_3)$ reduces to a 2-term pre-Lie$_\infty$ algebra or equivalently a pre-Lie $2$-algebra.
\end{rmk}

In \cite{Sheng}, the author showed that  skeletal pre-Lie $2$-algebras are classified by the third cohomology group of pre-Lie algebras. See~\cite{BC} also for more details for the classification of skeletal Lie 2-algebras. Similarly, we have

\begin{thm}
  There is a one-to-one correspondence between $2$-term skeletal \opt homotopy post-Lie algebras and triples $((\h,[\cdot,\cdot]_\h,\rhd),(V;\rho,\mu,\nu),\omega)$, where $(\h,[\cdot,\cdot]_\h,\rhd)$ is a post-Lie algebra, $(V;\rho,\mu,\nu)$ is a representation of the post-Lie algebra $(\h,[\cdot,\cdot]_\h,\rhd)$ and $\omega\in\Hom(\wedge^2\h\otimes \h,V)$ is a $3$-cocycle of $(\h,[\cdot,\cdot]_\h,\rhd)$ with coefficients in $(V;\rho,\mu,\nu)$.
\end{thm}
\begin{proof}
  Let $(\g=\g_{0}\oplus \g_{-1},[\cdot,\cdot]_\g,\Oprn_1,\Oprn_2,\Oprn_3)$ be a 2-term skeletal \opt homotopy post-Lie algebra, i.e. $\Oprn_1=0$. Then by Condition {\rm$(e_1)$} in Definition \ref{defi:2-term}, we deduce that $(\g_0,[\cdot,\cdot]_\g,\Oprn_2)$ is a post-Lie algebra.  Define linear maps  $\rho,~\mu,~\nu$ from $\g_0$ to $\gl(\g_{-1})$   by
  $$
  \rho(x)(a):=[x,a]_\g,\quad \mu(x)(a):=\Oprn_2(x,a),\quad \nu(x)(a):=\Oprn_2(a,x),\quad\forall x\in\g_0,~a\in\g_{-1}.
  $$ Obviously, $(\g_{-1};\rho)$ is a representation of the Lie algebra $(\g_0,[\cdot,\cdot]_\g)$. Then by   {\rm$(b_2)$},  {\rm$(b_3)$}, {\rm$(e_2)$} and {\rm$(e_3)$} in Definition \ref{defi:2-term}, we deduce \eqref{rep-1}-\eqref{rep-4} respectively. Thus,  $(\g_{-1};\rho,\mu,\nu)$ is a representation of the post-Lie algebra $(\g_0,[\cdot,\cdot]_\g,\Oprn_2)$. Finally, by  {\rm$(c)$} and {\rm$(f)$} in Definition \ref{defi:2-term}, we deduce that $\Oprn_3$ is a 3-cocycle of the  post-Lie algebra $(\g_0,[\cdot,\cdot]_\g,\Oprn_2)$ with coefficients in the representation  $(\g_{-1};\rho,\mu,\nu)$.

The proof of the other direction is similar. So the details will be omitted. \end{proof}

\vspace{2mm}
\noindent
{\bf Acknowledgements. } This research is supported by NSFC (11471139, 11425104, 11771190) and NSF of Jilin Province (20170101050JC). C. Bai is also supported by the Fundamental Research Funds for the Central Universities and Nankai ZhiDe Foundation.

 \end{document}